\documentclass[12pt]{amsart}

\usepackage{amssymb}
\usepackage{amsmath}
\usepackage{amscd}
\usepackage{graphicx}
\usepackage{amsfonts}

\numberwithin{equation}{section}

\newtheorem{thm}{Theorem}
\newtheorem{lem}{Lemma}
\newtheorem{cor}{Corollary}
\newtheorem{prop}{Proposition}

\newtheorem{defn}{Definition}

\newtheorem{rem}{Remark}

\begin{document}

\title[$p$-adic framed braids]{$p$-adic framed braids}

\author{J. Juyumaya}
\address{Universidad de Valpara\'{\i}so \\
Gran Breta\~na 1091, Valpara\'{\i}so, Chile.}
\email{juyumaya@uv.cl}

\author{S. Lambropoulou}
\address{ Department of Mathematics,
National Technical University of Athens,
Zografou campus, GR-157 80 Athens, Greece.}
\email{sofia@math.ntua.gr}
\urladdr{http://www.math.ntua.gr/$\tilde{~}$sofia}

\thanks{The first author was partially supported by Fondecyt 1050302, Dipuv and the NTUA. The second 
author was partially supported by Fondecyt (International Cooperation), U. Valparaiso and the NTUA}

\keywords{Inverse limits, $p$-adic integers, $p$-adic framed braids, Yokonuma-Hecke algebras, Markov
traces, $p$-adic framed links}

\subjclass{20C08, 57M27}

\date{}
\maketitle

\begin{abstract}
In this paper we define the $p$-adic framed braid group ${\mathcal F}_{\infty ,n}$, arising as the
inverse limit of the modular framed braids, and we give topological generators for
${\mathcal F}_{\infty ,n}$. We also give geometric interpretations for the $p$-adic 
framed braids. We then construct a $p$-adic Yokonuma-Hecke algebra ${\rm Y}_{\infty,n}(u)$ as the
inverse limit of a family of classical Yokonuma-Hecke algebras. These are quotients of the modular
framed braid groups over a quadratic relation. We also give topological generators for ${\rm
Y}_{\infty,n}(u)$. Finally, we  construct on this new algebra a linear trace that supports
 the Markov property. Paper presented at the 1017 AMS meeting.
\end{abstract}

\vspace{3mm} 
\begin{center}
{\sc Introduction}
\end{center}

\subsection{} 

Framed knots and links are like classical knots and links but with an integer, the
`framing', attached to each component. It is well--known that framed links can be used for constructing
$3$--manifolds using a topological technique called surgery \cite{li, wa}. Then two manifolds are 
homeomorphic if and only if any two framed links in $S^3$ representing them are related through isotopy 
moves and the Kirby moves (or the equivalent Fenn-Rourke moves) \cite{ki, fr}. In \cite{ks} K.H. Ko and L.
Smolinsky give a Markov-type equivalence for framed braids corresponding to homeomorphism classes of
$3$--manifolds. It would be certainly very interesting if one could construct $3$--manifold invariants by
constructing Markov traces on quotient algebras of the framed braid group and using the framed braid
equivalence of \cite{ks}. 

\smallbreak

 In 2.1 we recall the structure  of the framed braid group ${\mathcal F}_{n}= {\Bbb Z}^n\rtimes B_n$, where
$B_n$ is the classical braid group on $n$--strands. By construction, a framed braid splits into the `framing
part' and the `braiding part'. Moreover, ${\mathcal F}_{n}$ is generated by the elementary braids $\sigma_1,
\ldots, \sigma_{n-1}$ and by the elementary framings $f_1, \ldots, f_n$. We further introduce the modular
framed braid group  ${\mathcal F}_{d,n} = ({\Bbb Z}/d{\Bbb Z})^n\rtimes B_n$, which  has the same
presentation as ${\mathcal F}_n$, but with the additional relations:
$$
f_i^d = 1.
$$

In \cite{yo} the Yokonuma--Hecke algebras (shortened as Y--H algebras), ${\rm Y}_{d,n}(u)$, were introduced
by T. Yokonuma, where $u$ is a fixed non-zero complex number. They appeared originally in the representation
theory of finite Chevalley groups and they are natural generalizations of the classical Iwahori-Hecke
algebras, see also \cite{th}. In \cite{ju} Markov traces have been constructed by the first author for the
Y--H algebras of any index.  In Section 3 we introduce the  Y--H  algebra as a finite dimensional
quotient of the  group algebra ${\Bbb C}{\mathcal F}_{d,n} $ of ${\mathcal F}_{d,n}$ over the quadratic
relations:
$$
g_i^2 = 1 + (1-u)e_{d,i} (1- g_i), 
$$
where $g_i$ is the generator associated to the elementary braid $\sigma_i$ and  $e_{d,i}$ are certain
idempotents in ${\Bbb C}{\mathcal F}_{d,n}$, see 3.1. In ${\rm Y}_{d,n}(u)$ the relations
$f_i^d = 1$ still hold, and they are essential for the existence of the idempotents 
$e_{d,i}$, because $e_{d,i}$ is by definition a sum involving all powers of $f_i$ and $f_{i+1}$. In 3.4 we
give diagrammatic interpretations for the elements $e_{d,i}$ as well as for the quadratic relation
(see Figures~10, 11 and 12). 

\smallbreak

For relating to framed links and $3$--manifolds we would rather not have the restrictions $f_i^d = 1$ on
the framings. An obvious idea would be to consider the quotient of the classical framed braid group
algebra, ${\Bbb C}{\mathcal F}_{n}$, over the above quadratic relations. But then, the elements $e_{d,i}$
are not well-defined.
 Yet, we achieve this aim by employing the construction of inverse limits. In 1.1 we give some preliminaries
on inverse systems and inverse limits and we introduce the concept of topological generators. This is a set,
whose span is dense in the inverse limit (see Definition 1). In the rest of Section 1 we focus on
the construction of the $p$-adic integers ${\Bbb Z}_p$ and their approximations. 

\bigbreak

We shall now explain briefly our constructions.  Let $p$ be a prime number and  let $C_r$ be the cyclic 
group  of $p^r$ elements: $C_r\cong {\Bbb Z}/p^r{\Bbb Z}$. Then $\varprojlim C_r = {\Bbb Z}_p$, where the
inverse system maps  $\theta_s^r : {\Bbb Z} / p^r{\Bbb Z} \longrightarrow
 {\Bbb Z} / p^s{\Bbb Z}$  ($r\geq s$) are the natural epimorphisms.
${\Bbb Z}_p$ contains ${\Bbb Z} = \langle{\bold t} \rangle$ as a dense subgroup. The element ${\bold t}$ is
a topological generator for ${\Bbb Z}_p$ (see subsection 2.2), and a $p$-adic integer will be denoted 
$ {\bold t}^{\underleftarrow{a}}$, where $\underleftarrow{a} =: (a_1, a_2,\ldots)$ with $a_r \equiv a_s \
(\text{mod} \, p^s)$ whenever $r\geq s$.  We consider the inverse system  $\left( C_r^n, \pi^r_s \right)$
indexed by ${\Bbb N}$, where the map
$\pi_r^s: C_r^n \longrightarrow C_s^n$ ($r\geq s$) acts componentwise as the natural epimorphism
$\theta_s^r$. Then
$\varprojlim C_r^n \cong {\Bbb Z}_p^n$ (see Proposition~\ref{ciarn}) and ${\Bbb Z}_p^n$ contains 
${\Bbb Z}^n = \langle {\bold t}_1, \ldots, {\bold t}_n \rangle$ as a dense subgroup (Lemma~\ref{gencrn}). 
 We then consider the inverse system $\left( {\mathcal F}_{p^r,n}, \, \pi^r_s\cdot {\rm id}\right)$ indexed
by ${\Bbb N}$, where the map $\pi^r_s\cdot {\rm id}$ acts on the framing part of a modular
framed braid as described above, and trivially on the braiding part. So, we define the {\it $p$--adic
framed braid group} ${\mathcal F}_{\infty , n}$ (see Definition \ref{finfty}) as 
$$
{\mathcal F}_{\infty , n} =  \varprojlim {\mathcal F}_{p^r,n}.
$$
All this is explained in Section 2, subsections 2.2, 2.3 and 2.4. Geometrically, a $p$-adic framed braid is
an infinite sequence of modular framed braids with the same braiding part and such that the framings of the
$i$th strand in each element of the sequence give rise to a $p$-adic integer. See subsection 2.5 and 
right-hand side of Figure~5 for an illustration. So, a $p$-adic framed braid can be also interpreted as an
infinite framed cabling of a braid in $B_n$, such that the framings of each infinite cable form a $p$-adic
integer. See Figure~1, where $ (a_1, a_2,\ldots), (b_1, b_2,\ldots)  \in {\Bbb Z}_p$.

\smallbreak 

\begin{figure}[h]
\begin{center}
\includegraphics[width=5in]{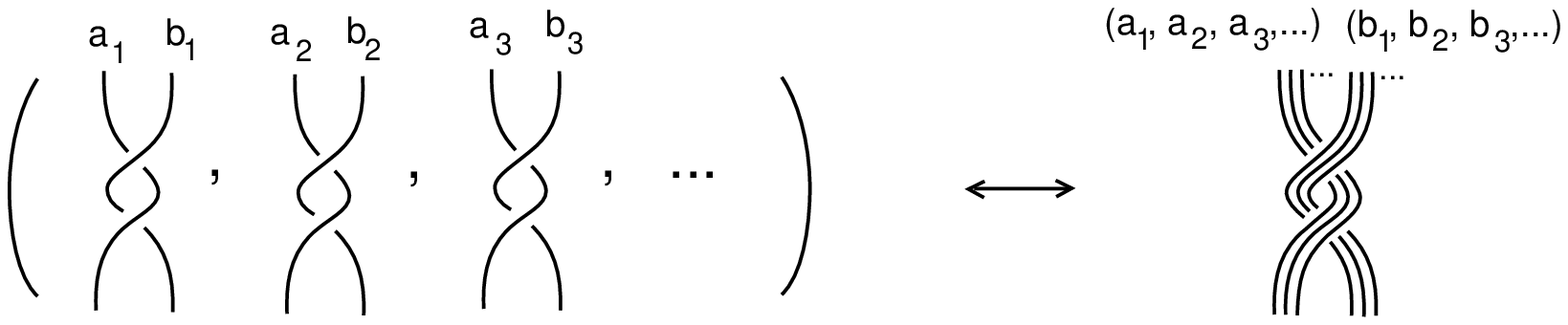}
\end{center}
\caption{A $p$-adic framed braid as an infinite framed cabling}
\label{figure1}
\end{figure}

In Theorem~1 the natural identification
$$
{\mathcal F}_{\infty , n} \cong {\Bbb Z}_p^n\rtimes B_n
$$ 
is established. This identification implies, in particular,
that there are no modular relations for the framing in ${\mathcal F}_{\infty , n}$. Also,
that the classical framed braid group ${\mathcal F}_n$ sits in ${\mathcal F}_{\infty , n}$ as a dense
subset. Hence,  the set $A = \{ {\bold t}_1,\sigma_1,\ldots, \sigma_{n-1} \} \subset {\mathcal F}_n$
 is a set of topological generators for ${\mathcal F}_{\infty , n}$ (see Theorem~1). The 
identification  ${\mathcal F}_{\infty , n} \cong {\Bbb Z}_p^n\rtimes B_n$ is also used in 2.5, where we give
geometric interpretations of the $p$-adic framed braids as classical braids with framings $p$-adic
integers. See Figure~5. We can say that a $p$-adic framed braid splits into the
`$p$-adic framing part' and the `braiding part'. That is, a $p$-adic framed braid is a word of the form:
$$
{\bold t}_1^{\underleftarrow{a_1}} {\bold t}_2^{\underleftarrow{a_2}} \cdots  {\bold
t}_n^{\underleftarrow{a_n}} \cdot \sigma
$$
where $\underleftarrow{a_1},\ldots, \underleftarrow{a_n}$ are the $p$-adic framings and $\sigma \in B_n$.
  Of course, the closure of a $p$-adic framed braid defines an oriented
$p$-adic framed link, see Figure~2 for an example.

\smallbreak 

\begin{figure}[h]
\begin{center}
\includegraphics[width=2.6in]{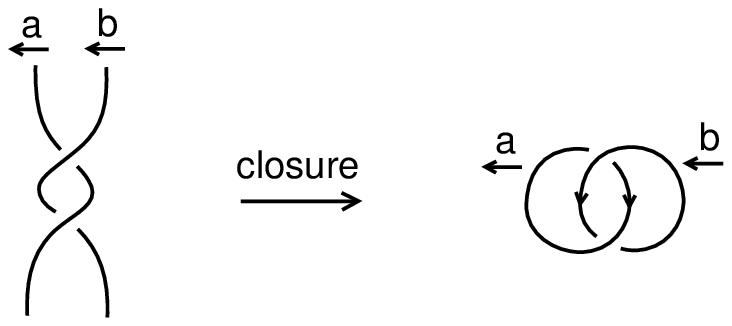}
\end{center}
\caption{A $p$-adic framed braid and a $p$-adic framed link}
\label{figure2}
\end{figure}

\noindent In 2.6 we give approximations of $p$-adic framed braids by sequences of classical
framed braids. See Figure~9 and Figure~10 for examples.

\smallbreak

 In  3.5 we define the {\it $p$--adic Yokonuma--Hecke algebra}  ${\rm Y}_{\infty,n}(u)$ 
as the inverse limit of the inverse system $\left({\rm Y}_{p^r,n}(u), \varphi_s^r\right)$ of classical Y--H
algebras, indexed by ${\Bbb N}$ (see Definition~\ref{yinfty}): 
$$
{\rm Y}_{\infty ,n}(u) = \varprojlim {\rm Y}_{p^r,n}(u).
$$
The above inverse system is induced by  the inverse system $\left( {\Bbb C}{\mathcal F}_{p^r,n},
\phi_s^r\right)$, where $\phi_s^r$ is the `linear span' of $\pi^r_s\cdot {\rm id}$ at the level of the
group algebra, using also our definition of the  Y--H  algebras as finite dimensional 
quotients of the  group algebras ${\Bbb C}{\mathcal F}_{d,n}$.  ${\rm Y}_{\infty,n}(u)$ is an infinite
dimensional algebra, in which the framing restrictions $f_i^d = 1$ do not hold.
 In 3.6, Theorem 3, we give a set of topological generators $\{ {\bf t_1}, g_1, \ldots ,g_{n-1}\}$ for ${\rm
Y}_{\infty,n}(u)$  satisfying the quadratic relations:
$$
g_i^2 = 1 + (1-u)e_i (1- g_i),
$$
 where the element $e_i$  is also an idempotent and its approximation involves the `framing' generators
${\bf t_i, t_{i+1}}$. 

\smallbreak

In Section 4 we recall the Markov traces on $\bigcup_{n=1}^{\infty}{\rm Y}_{d,n}$, constructed in
\cite{ju}, and using these traces we extend the construction to a {\it $p$-adic Markov trace} $\tau$ on
$\bigcup_{n=1}^{\infty}{\rm Y}_{\infty,n}$ to some inverse limit of polynomial rings. See Theorem~5. 
 Finally, in 4.3 we give some computations of the trace $\tau$ on identity braids with $p$-adic framings,
on the elements $e_i$ and  on $g_i^2$.

\smallbreak

It is, perhaps, worth noting that, adapting to the case of framed braids the cyclotomic and `generalized' 
Hecke algebras of  type B, used in \cite{la} by the second author for constructing Markov traces and link
invariants for the solid torus,  gives rise to some `framed' Hecke algebras, which are special cases of the
Y--H algebras. In all the above algebras the quadratic relation is the same as in the classical
Iwahori--Hecke algebra, so  in \cite{la} it was possible to remove the modular relation for the framing
and construct the generalized B--type Hecke algebras, without employing the theory of inverse limits.

\smallbreak

In a sequel paper we combine the $p$-adic trace $\tau$, as well as the Markov traces constructed in
\cite{ju},  with  the Markov  equivalence for $p$-adic framed braids to construct invariants of oriented
$p$-adic framed links. We hope that this new concept of  $p$-adic framed braids and  $p$-adic framed links
that we propose, as well as our $p$-adic framing invariant will be useful for constructing $3$--manifold
invariants using the theory of braids and the Markov-type equivalence given in \cite{ks}.

\subsection{}

As usual we denote by ${\Bbb C}$, ${\Bbb Z}$ and ${\Bbb N}=\{1, 2, \ldots \}$ the set of complex numbers,  
the integers and the natural numbers respectively. We also  denote  ${\Bbb Z} / d{\Bbb Z}$  the additive
group of integers modulo $d$.  Throughout the paper we fix a prime number $p$ and a $u$ in ${\Bbb
C}\backslash\{0\}$. Finally, whenever two elements $a, b$ are identified we shall write
$a \stackrel{!}{=} b$.

\subsection{}

 Let $H$ be a group and let $H^n= H\times \cdots \times H $ ($n$--times). The symmetric group $S_n$ of the
permutations of the set $\{1,2, \ldots, n\}$ acts on $H^n$ by permutation, that is:
$$
\sigma \cdot (h_1,\ldots,h_n) = (h_{\sigma(1)},\ldots,h_{\sigma(n)}) \ \ \ \ \forall \ \sigma \in S_n.
$$
We define on the set $H^n \times S_n$ the operation: 
$$
(h, \sigma) \cdot (h^{\prime}, \tau) = (h \, \sigma (h^{\prime}), \sigma \, \tau).
$$
Then, the set $H^n \times S_n$  with the above operation is a group, called {\it the wreath product} of $H$
and $S_n$, denoted by $H  \wr S_n$. Recall that $H  \wr S_n = H^n \rtimes  S_n$.

\section{Inverse Limits  and p-adic Integers}

\subsection{}

 An {\it inverse system} $(X_i, \phi_j^i)$ of topological spaces indexed by a directed set $I$, consists of
a family $(X_i \ ; \ i\in I)$ of topological spaces and a family 
$(\phi_j^i: X_i \longrightarrow X_j \ ; \  i,j \in I, \ i\geq j)$ of continuous maps, such that 
$$ 
\phi_i^i = {\rm id}_{X_i} \quad{\rm and} \quad  \phi_k^j \circ \phi_j^i  = \phi_k^i 
\quad {\rm whenever}  \quad
i \geq j \geq k.
$$
Sets with no other topology specified are regarded as topological spaces  with the discrete
topology. In particular, finite sets are compact Hausdorff spaces. If each $X_i$ is a topological group
the continuous maps
$\phi_j^i$ are required to be also homomorphisms. Inverse systems of topological rings, topological vector
spaces, topological algebras, et cetera are defined analogously.

\smallbreak
The {\it inverse limit} $\varprojlim X_i$ of the inverse  system $(X_i, \phi_j^i)$  is defined to be  the
following  subset of the cartesian product $\prod X_i$.
$$
\varprojlim X_i := \{ z \in \prod X_i\, ; \,  (\phi_j^i\circ \varpi_i )(z) = \varpi_j (z)
\quad \text{ whenever}  \quad j\geq i\},
$$
where the map $\varpi_i$ denotes the natural projection of
$\prod X_i$ onto $X_i$. It turns out that $\varprojlim X_i$ is uniquely defined, and it is non-empty if
each $X_i$ is a non-empty compact Hausdorff space.
 If the $X_i$'s are topological groups then $\varprojlim X_i$ is a topological group with operation
induced in $\prod X_i$ componentwise by the group operations. Moreover, in this case, $\varprojlim X_i$ is
always non-empty.  Similar facts are true for  topological rings, topological vector spaces, topological
algebras, et cetera. We also note that, if $X_i =X$ for all $i$ and $\phi_j^i$ is the identity for all $i,j$ then
$\varprojlim X$ can be identified naturally with $X$ (identifying a constant sequence $(x,x,\ldots)$ with $x\in
X$).

\smallbreak
As a topological space, $\prod X_i$ is endowed with the product topology, so $\varprojlim X_i$ inherits the
induced topology. It can be then verified  that $\varprojlim X_i$ is closed in  $\prod X_i$. A {\it basis}
of open sets in $\varprojlim X_i$  contains elements of the form 
$$
  \varpi_{i}^{-1}(U_i)  \cap \varprojlim X_i,
$$
where  $U_i$ open in $X_i$. Then, any open set in $\varprojlim X_i$ is a union of sets of the form 
 \begin{equation} \label{open} 
\varpi_{i_1}^{-1}(U_1) \cap \cdots \cap \varpi_{i_n}^{-1}(U_n)  \cap \varprojlim X_i,
 \end{equation}
where $i_1, \ldots,i_n
\in I$ and $U_r$ open in $X_{i_r}$ for each $r$. (Compare with \cite{riza}, p.7.) 

\smallbreak
A {\it morphism }  between two  inverse systems $(X_i, \phi_j^i)$ and $(Y_i, \psi_j^i)$, both   indexed by
the same directed set $I$, is a collection of continuous maps 
$$
(\rho_i: X_i \longrightarrow Y_i \ ; \ i\in I)
$$  
such that $\psi_j^i \circ \rho_i = \rho_j \circ \phi_j^i $, for all $i\in I$. If, moreover, the inverse
systems are topological groups then the maps $\rho_i$ must be also group homomorphisms. The definition
is analogous if the inverse systems considered are topological rings, et cetera.
 A morphism/group homomorphism   $(\rho_i \ ; \ i\in I)$  from the inverse system  $(X_i, \phi_j^i)$ to
the inverse system  $(Y_i, \psi_j^i)$ induces a morphism/group homomorphism between the inverse limits:
$$
\varprojlim \rho_i :  \varprojlim X_i \longrightarrow \varprojlim Y_i
$$ 
by setting
$$
\varprojlim \rho_i ((x_i)) := (\rho_i (x_i)).
$$
If we have embeddings $\iota_i$ from $ X_i$ into  $Y_i$, these induce  a natural
embedding  $\varprojlim \iota_r  : \varprojlim X_i  \longrightarrow \varprojlim Y_i$.
Moreover, if the following sequence    
$$
\begin{CD}
0 @>>>X_i @>\iota_i>> Y_i @>\varphi_i>>Z_i@>>> 0
\end{CD}
$$ 
is exact for any $i$, then  the sequence 
\begin{equation}\label{exact}
\begin{CD}
0 @>>>\varprojlim X_i  @> \varprojlim \iota_r>>  \varprojlim Y_i@>\varprojlim\varphi_r>>
\varprojlim Z_i
\end{CD}
\end{equation}
is  exact (see \cite{bou, riza}). 
\smallbreak

Let now $J$ be a {\it cofinal subset}  of the index set $I$ (that is, for every $i \in I$ there is a $j \in
J$ with $j\geq i$). Then $J$ gives rise to the same inverse limit, that is, $\varprojlim_{i\in I}X_i
\cong \varprojlim_{j\in J}X_j$. Finally, let  $X$ and $Y$ be the  inverse limits of the inverse systems $(X_i,
\phi_k^i \ ; \ i\in I)$ and $(Y_j, \psi_m^j \ ; \ j\in I)$, respectively.  Then we have 
\begin{equation}\label{limpro}
X\times Y \cong \varprojlim_{(i,i)}(X_i\times Y_i) \cong \varprojlim_{(i,j)\in I\times I}(X_i\times Y_j).
\end{equation}
 The bijection between $X\times Y $ and
$\varprojlim_{(i,i)}(X_i\times Y_i)$ identifies pairs of sequences
$((x_i), (y_i)) \in  X\times Y$ with the sequence $(x_i, y_i) \in  \varprojlim_{(i,i)}(X_i\times Y_i)$, and
it becomes an isomorphism when  $X_i, Y_i$ are groups, rings, et cetera. Note here that the diagonal
subset is a cofinal subset of $I\times I$.  Clearly, the above generalize to any finite cartesian product
of inverse limits.

\subsection{}

Our working example for the notion of inverse limit will be the construction of the $p$-adic integers.
 Let $p$ be a prime number, which will be fixed throughout the paper, and let ${\Bbb Z} / p^i{\Bbb Z}$ be
the additive group of integers modulo $p^i$.  For any $r, s \in {\Bbb N}$ with $r\geq s$ we
consider the following natural epimorphisms:
\begin{equation}\label{thetarsadd}
\begin{array}{cccc}
\theta_s^r : &  {\Bbb Z} / p^r{\Bbb Z} &  \longrightarrow & {\Bbb Z} / p^s{\Bbb Z} \\
            \,&  k + p^r{\Bbb Z} & \mapsto &   k + p^s{\Bbb  Z}
\end{array}
\end{equation}
More precisely, an element $a_r \in {\Bbb Z} / p^r{\Bbb Z}$ can be written uniquely in the form
$$
a_r = k_0 + k_1 p + k_2 p^2 + \cdots + k_{r-1} p^{r-1} + p^r{\Bbb Z}, \ \text{where} \, \
k_0,\ldots,k_{r-1} \in
\{0,1,\ldots,p-1\}. 
$$
Then, for $r\geq s$: 
$$
\theta_s^r(k_0 + k_1 p + k_2 p^2 + \cdots + k_{r-1} p^{r-1} + p^r{\Bbb Z}) =  k_0 + k_1 p + k_2 p^2 + \cdots +
k_{s-1} p^{s-1} + p^s{\Bbb Z}
$$
(``cutting out'' $r-s$ terms). 
 We obtain, thus, the inverse system  $({\Bbb Z} / p^r{\Bbb Z} , \theta_s^r)$ of topological groups, indexed
by ${\Bbb N}.$ Its inverse limit, $\varprojlim {\Bbb Z} / p^r{\Bbb Z}$, is  the group of
{\it $p$-adic integers}, denoted  ${\Bbb Z}_p$. 
 Clearly, an entry $a_r \in {\Bbb Z} / p^r{\Bbb Z}$ of an element $(a_r) \in {\Bbb Z}_p$ can be identified with
\begin{equation}
\label{padic}
a_r = k_0 + k_1 p + k_2 p^2 + \cdots + k_{r-1} p^{r-1},
\end{equation}
where
$k_0,\ldots,k_{r-1} \in \{0,1,\ldots,p-1\}$. Thus, ${\Bbb Z}_p$ can be identified with the set of power series:
\begin{equation}\label{series}
{\Bbb Z}_p = \{ \sum_{n=0}^{\infty} b_n p^n   \ ; \ b_n \in  {\Bbb N}, \ 0\leq b_n < p\}.
\end{equation}
Equivalently, 
${\Bbb Z}_p$ can be identified with the set of sequences:
\begin{equation}\label{zipi}
{\Bbb Z}_p = \{ (a_r) \ ; \ a_r \in  {\Bbb Z}, \ a_r \equiv a_s \ (\text{mod} \, p^s) \ \ \text{whenever} \
r\geq s\}.
\end{equation}
Elements in ${\Bbb Z}_p$ shall be usually denoted as 
\begin{equation}\label{parrow}
\underleftarrow{a} := (a_1, a_2, a_3,\ldots) \in {\Bbb Z}_p
\end{equation}
${\Bbb Z}_p$ is a  subgroup of $\prod {\Bbb Z} / p^r{\Bbb Z}$ with operation
inherited  by the componentwise (additive) operation of $\prod {\Bbb Z} / p^r{\Bbb Z}$. ${\Bbb Z}_p$ is  non-cyclic
and it contains no elements of finite order.  Each ${\Bbb Z} / p^r{\Bbb Z}$ is  finite, making  ${\Bbb Z}_p$  {\it
profinite} (see \cite{riza}, \cite{wi}).

\subsection{}

Contrary to embeddings between inverse systems, if each component  $\rho_i: X_i \longrightarrow Y_i$ of a
morphism between two inverse systems is onto, the induced map  $\varprojlim \rho_i$ between the inverse
limits is not necessarily onto.  

\smallbreak
For example, consider the inverse systems  $({\Bbb Z}, \text{id})$
and  $({\Bbb Z} / p^r{\Bbb Z} , \theta_s^r)$, both indexed by ${\Bbb N}$, and for each $s \in {\Bbb N}$
define  the canonical  epimorphism
\begin{equation}\label{canepi}
\begin{array}{cccc}
\rho_s: &  {\Bbb Z} &  \longrightarrow & {\Bbb Z} / p^s{\Bbb Z} 
\end{array}
\end{equation}
To see $\rho_s$ concretely, let $a\in {\Bbb Z}$. Then $a$ can be written uniquely in the form
$a = k_0 + k_1 p + \cdots + k_i p^i$, for some $i \in {\Bbb N}$, where $k_0,\ldots,k_i \in
\{0,1,\ldots,p-1\}$. So, 
$$
\begin{array}{cclc}
\rho_s(a) &  = &   k_0 + k_1 p + \cdots + k_i p^i & \text{for} \ i < r \\
\rho_s(a) &  = &   k_0 + k_1 p + \cdots + k_{s-1} p^{s-1}  = \theta_s^{i+1}(a +  p^{i+1}{\Bbb Z}) & \text{for} \ 
i \geq s
\end{array}
$$
 Then $(\rho_s \ ; \ s \in {\Bbb N})$ is a homomorphism between the two inverse systems.
The first inverse limit is isomorphic to ${\Bbb Z}$ (identify $(z,z,\ldots) \in \varprojlim {\Bbb Z}$ with
$z \in {\Bbb Z}$), while the second is the set of 
$p$-adic integers ${\Bbb Z}_p$.  Note that the image of  $\varprojlim {\Bbb Z}$ in ${\Bbb Z}_p$ under
$\varprojlim\rho_s$ consists in all constant tuples  of integers (or, according to the notation of (\ref{series}),
it consists in the finite power series). On the other hand, the tuple 
 $(b_r)$,  such that $b_r = 1+p+ \cdots + p^{r-1} $ is in  ${\Bbb Z}_p$ but is not constant. 

\smallbreak
Yet, we have the following very important result.

\begin{lem}{\rm(Lemma 1.1.7. in \cite{riza})}  Let $(X_i, \phi_j^i)$ be an inverse system of
topological spaces indexed by a directed set $I$ and let $\rho_i :   X \longrightarrow  X_i$ be compatible
surjections from a topological space $X$ onto the spaces $X_i$ ($i \in I$). Then, either $\varprojlim X_i =
\emptyset$ or the induced mapping  $\rho = \varprojlim\rho_i :   \varprojlim X \longrightarrow  \varprojlim X_i$
maps $\varprojlim X$ onto a dense subset of $\varprojlim X_i$.
\end{lem}
\begin{proof} 
For the proof of Lemma 1 consider a non--empty open set $V$ in $\varprojlim X_i$ of the form
(\ref{open}).  We have to show that $\rho(X) \cap V \neq \emptyset$. Indeed, let $i_0 \geq
i_1,\ldots,i_n$ and let $y= (y_i) \in V$. Choose $x\in X$ so that $\rho_{i_0} (x) = y_{i_0}$. Then $\rho
(x)  \in V$. 
\end{proof}

For example, let $\rho_i$ denote the restriction  of the canonical projection of $\varprojlim X_i$ onto
$X_i$ on  a subset $A \subset \varprojlim X_i$.  Recall that $\varprojlim A$ can be identified with $A$. Then we
have the following.
\begin{cor}\label{dense1}
If for a subset $A \subset \varprojlim X_i$ we have $\rho_i(A) = X_i$ for all $i\in I$,  then $\rho (\varprojlim
 A)$  is dense in  $\varprojlim X_i$, where $\rho = \varprojlim\rho_i$.
\end{cor}
Since ${\Bbb Z}$ projects onto each factor ${\Bbb Z} / p^r{\Bbb Z}$ via the canonical epimorphism (\ref{canepi}),
we obtain the following, as an application of Corollary \ref{dense1}.

\begin{cor}\label{dense2}
  ${\Bbb Z}$ is dense in ${\Bbb Z}_p$.
\end{cor}
This means that every $p$-adic integer can be approximated by a sequence of constant sequences. In 1.4 we study
further this approximation.
\begin{defn}\label{generators}
\rm (cf. \S \, 2.4 in \cite{riza}) Let $G_i$ be a group (ring, algebra, et cetera) for all $i\in I$. A subset $S
\subset
\varprojlim G_i$  is a set of {\it topological generators} of $\varprojlim G_i$ if the span $\langle S\rangle$ is
dense in 
$\varprojlim G_i$. If, moreover, 
$S$ is finite, $\varprojlim G_i$ is  said to be {\it finitely generated}.
\end{defn}
For example, the element $(1,1,\ldots)$ is a topological generator of ${\Bbb Z}_p$, since, by Corollary
\ref{dense2}, the cyclic subgroup $\langle (1,1,\ldots) \rangle = \Bbb Z$ is dense in ${\Bbb Z}_p$.

\subsection{}

Let us now take a closer look at the $p$-adic integers. Seen as a tuple of integers, an element of ${\Bbb Z}_p$
may have many disguises. For example: 
$$
(1,\, p+1,\, p^2+p+1,\, a_4, a_5, \ldots) = (p^2+p+1,\, p^2+p+1,\, p^2+p+1,\, a_4, a_5, \ldots),  
$$
since $ p+1 \equiv p^2+p+1 \, (\text{mod} \, p^2)$ and $ 1 \equiv p^2+p+1 \, (\text{mod} \, p)$. From (\ref{zipi})
and from the form (\ref{padic}) of an entry $a_r \in {\Bbb Z} / p^r{\Bbb Z}$ of an element $(a_r) \in {\Bbb Z}_p$,
it follows that for the  ($n+1$)st entry there are $p$ choices, namely: 
\begin{equation}
\label{ar+1}
a_{r+1} \in \{a_r + \lambda p^r \ ; \ \lambda= 0,1,\ldots, p-1 \}.
\end{equation}
 On the contrary, there is no choice for the entries before, as $a_s \equiv a_r (\text{mod}\, p^s)$ for all 
$s=1,\ldots,n-1$.

\smallbreak
On the level of basic open sets the logic is similar: As a topological space, ${\Bbb Z}_p$ is endowed with the
induced topology of $\prod ({\Bbb Z} / p^r{\Bbb Z})$, which builds up from the discrete topology of each factor 
${\Bbb Z} / p^r{\Bbb Z}$. Thus, a basic open set in ${\Bbb Z}_p$ is of the form 
 $\{ \varpi_i^{-1} (U_i) \ ; \ U_i \subseteq {\Bbb Z} / p^i{\Bbb Z} \}$, where $\varpi_i$
is the restriction of the natural projection of ${\Bbb Z}_p$ onto $ {\Bbb Z} / p^i{\Bbb Z}$. So,  we have for
example:
\[ 
\begin{array}{lcl}
\varpi_1^{-1}(\{1\}) & = &   (\{1\}\times {\Bbb Z} / p^2{\Bbb Z} \times {\Bbb Z} / p^3{\Bbb Z} \times \cdots )  
                               \cap {\Bbb Z}_p \\
             &  = &   \{ (1,a_2, a_3,a_4, \ldots) \in {\Bbb Z}_p \} \\
             &    &    \\
\varpi_2^{-1}(\{p+1\}) & = &  ({\Bbb Z} / p {\Bbb Z}\times \{p+1\} \times {\Bbb
Z} / p^3{\Bbb Z}\times  \cdots )   \cap   {\Bbb Z}_p \\
             &  = &   \{ (1,\, p+1,\, a_3, a_4, \ldots) \in {\Bbb Z}_p \}  \\
             &    &    \\
\varpi_3^{-1}(\{p^2+p+1\})  &  =  & \{ (1,\, p+1,\, p^2+p+1,\, a_4, \ldots) \in {\Bbb Z}_p \}, \\
\end{array}
\] 
where the $a_i$'s are subject to the conditions (\ref{ar+1}). For $U_i$  not a singleton, $\varpi_i^{-1} (U_i)$ is
the union of the sets $\varpi_i^{-1} (\{u\})$, for all $u \in U_i$. In the above example note that 
$$
\varpi_3^{-1}(\{ p^2+p+1\}) \subsetneq \varpi_2^{-1}(\{ p+1 \}) \subsetneq \varpi_1^{-1}(\{ 1 \}).
$$
 In general, for any given element  $a= (a_1,a_2,\ldots) \in {\Bbb Z}_p$ we have that 
$a \in \varpi_i^{-1}(\{a_i\})$ for all $i$. Moreover,
\[ 
\begin{array}{l}
\varpi_1^{-1}(\{a_1\})  =   \{ (a_1,x_2, x_3, x_4, \ldots) \ ; \ x_2 \equiv a_1 \, (\text{mod} p), \ x_n \equiv
x_m \, (\text{mod} p^m), \ n\geq m \} \\
                     \\
\varpi_2^{-1}(\{a_2\})  =   \{ (a_1,a_2, y_3, y_4, \ldots) \ ; \ y_2 \equiv a_2 \, (\text{mod} p^2), \ y_n
\equiv y_m \, (\text{mod} p^m), \ n\geq m \} \\
 \ \ \ \ \ \ \ \ \ \ \ \ \ \ \ \vdots      \\
\end{array}
\] 
As we can see, $ \varpi_1^{-1}(\{ a_1 \}) \supsetneq \varpi_2^{-1}(\{a_2 \}) \supsetneq \cdots.$

\bigbreak
This implies, in particular, that for any $k \in  {\Bbb Z}$ the constant sequence $(k,k,\ldots) \in  {\Bbb Z}_p$
 is contained in infinitely many open sets, each set being a refinement of the previous one. Recall now, from
Corollary~\ref{dense2}, that the set of constant sequences is dense in ${\Bbb Z}_p$. From the above, for every
element $\underleftarrow{a} \in {\Bbb Z}_p$ we can construct an appproximating sequence of constant sequences
in ${\Bbb Z}_p$. 
 For example, the  sequence referred to in 1.3:
\begin{equation}
\label{bn}
\underleftarrow{b} = (1, \, 1+p, \, 1+p+p^2, \, \ldots)
\end{equation}
 can be approximated in ${\Bbb Z}_p$ by the following sequence of constant sequences:
\[ 
\begin{array}{cl}

(1, \, 1, \, 1, \, \ldots) & \in \varpi_1^{-1}(\{1\}) \\
               &   \\
(1+p,\, 1+p,\, 1+p,\, \ldots) & \in \varpi_2^{-1}(\{1+p\})   \\
               & \\
(1+p+p^2,\, 1+p+p^2,\, 1+p+p^2,\, \ldots) \ &  \in \varpi_3^{-1}(\{1+p+p^2\})  \\
           \vdots  & \\
\end{array} 
\]
Indeed, $(b_n) \in \varpi_i^{-1}(\{b_i\})$ for all $i$. Also, from the observation in the beginning of the
subsection, the above constant sequences have the same first entry, equal to $b_1$. Likewise, all  sequences from
the second one  have the same second entry, equal to $b_2$,  and so on. Finally, $ \varpi_1^{-1}(\{ 1 \})
\supsetneq \varpi_2^{-1}(\{1+p \}) \supsetneq \cdots.$

\smallbreak
In general, the element
$\underleftarrow{a} = (a_1, a_2, a_3, \ldots) \in {\Bbb Z}_p$ can be approximated by the following sequence of
constant sequences:
\begin{equation}
\label{appx}
\begin{array}{ll}
(a_1, a_1, a_1, \ldots) & \in \varpi_1^{-1}(\{a_1 \}) \\
               &   \\
 (a_1, a_2, a_2, a_2, \ldots) = (a_2, a_2, a_2, \ldots) & \in \varpi_2^{-1}(\{a_2 \}) \\
               & \\
 (a_1, a_2, a_3, a_3, \ldots) = (a_3, a_3, a_3, \ldots)  & \in \varpi_3^{-1}(\{a_3 \}) \\
 & \\
 \, \ \ \ \ \ \ \ \  \vdots  & \\
\end{array} 
\end{equation}
since  $a_r \equiv a_1 (\text{mod}\,\, p),$ for $r\geq 1$,  $a_r \equiv a_2 (\text{mod}\,\, p^2)$, for $r\geq 2$,
and so on.  Indeed,
$(a_1, a_2, a_3, \ldots) $ and $(a_i, a_i, a_i,
\ldots) =(a_1, a_2, \ldots, a_{i-1}, a_i, a_i, \ldots)$ are both in $\varpi_i^{-1}(\{a_i\})$ for all $i$.
Finally,
$\varpi_1^{-1}(\{a_1\}) \supsetneq \varpi_2^{-1}(\{ a_2 \})
\supsetneq \cdots$, justifying the approximation claim. We shall write:
\begin{equation}
\label{lim}
\underleftarrow{a} = \lim_{k} (a_k).
\end{equation}
 For more details and further reading on inverse limits and the $p$-adic integers see, for example,
\cite{riza, ro, wi}.

\section{$p$-adic framed braids}

The aim of this section is to introduce the notion of $p$-adic framed braids. These are similar to the
classical framed braids but, instead of integral framing, each strand may be coloured with any $p$-adic integer.

\subsection{}

Before starting with our construction we need to digress briefly and 
recall the definition and the structure of the classical framed braid group and the modular framed braid
group. 
\smallbreak
We consider the group ${\Bbb Z}^n$  with the usual operation:
\begin{equation}\label{additive}
(a_1, \ldots,a_n)(b_1, \ldots, b_n):= (a_1 + b_1, \ldots,a_n +b_n).
\end{equation}
${\Bbb Z}^n$ is generated by the `elementary framings': 
$$
f_i:= (0, \ldots,0,1,0, \ldots,0)
$$
with  $1$ in the $i$th position. Then, for example, $f_i^m = (0, \ldots,0,m,0, \ldots,0)$ and $f_i f_j = (0,
\ldots,0,1,0,
\ldots,0,1,0, \ldots,0),$ with  $1$ in the $i, j$ positions, and an element
$a=(a_1, \ldots,a_n)\in {\Bbb Z}^n$ can be expressed as:
$$
a= f_1^{a_1}f_2^{a_2}\cdots f_{n}^{a_n}.
$$

 Let also $B_n$ be the classical braid group on $n$ strands. $B_n$ is generated by the elementary braids $\sigma_1,
\ldots, \sigma_{n-1}$, where  $\sigma_i$ is the positive crossing between the $i$th and the
$(i+1)$st strand. The $\sigma_i$'s satisfy the well--known braid relations:
$\sigma_i\sigma_j  = \sigma_j\sigma_i$,  if  $ |i-j|>1 $ and 
$\sigma_i\sigma_{i+1}\sigma_i  = \sigma_{i+1}\sigma_i\sigma_{i+1}$. 
Recall the symmetric group $S_n$, generated by the $n-1$ elementary
transpositions $s_i := (i, i+1)$, and let further $\pi$ be  the natural projection of $B_n$
on $S_n$. We let  $\sigma(j)$ denote $\pi(\sigma)(j)$  for any $j = 1,2, \ldots, n$. In particular, 
$\sigma_i (j) = s_i(j).$ Using $\pi$ we define the {\it framed braid group} ${\mathcal F}_{n}$ (see 0.3) as: 
$$
{\mathcal F}_{n} = {\Bbb Z}  \wr B_n = {\Bbb Z}^n \rtimes  B_n,
$$
  where the action of $B_n$ on $a = (a_1 , \ldots,  a_n)\in {\Bbb Z}^n$ is given by permutation of the indices: 
\begin{equation}\label{action}
\sigma (a) =(a_{\sigma(1)} , \ldots, a_{\sigma(n)})\qquad (\sigma\in B_n).
\end{equation}
In the above notation, the action of $B_n$ on ${\Bbb Z}^n$ is given by the multiplicative formula:
$$
\sigma (f_1^{a_1}f_2^{a_2}\cdots f_{n}^{a_n})= f_{1}^{a_{\sigma(1)}}f_{2}^{a_{\sigma (2)}}\cdots
f_{n}^{a_{\sigma(n)}} \qquad (\sigma\in B_n).
$$
 
Any word in ${\mathcal F}_{n}$  splits, by construction, 
into the \lq framing\rq \, part and the \lq braiding\rq \, part. That is,  it can be written in the
form 
\begin{equation}\label{fsplit}
f_1^{k_1} f_2^{k_2}\cdots f_n^{k_n} \cdot \sigma,\quad \text{where} \quad k_i \in {\Bbb Z}, \ 
\sigma \in B_n.
\end{equation}
The multiplication in ${\mathcal F}_n$ is defined using the action of  $B_n$ on ${\Bbb Z}^n$ as follows:
\begin{equation}\label{frprod}
(f_1^{a_1}f_2^{a_2}\cdots f_{n}^{a_n}\cdot \sigma)(f_1^{b_1}f_2^{b_2}\cdots f_{n}^{b_n}\cdot \tau) :=
f_1^{a_1 + b_{\sigma(1)} }f_2^{a_2 + b_{\sigma(2)}}
\cdots f_{n}^{a_n + b_{\sigma(n)}} \cdot \sigma\tau .
\end{equation}
Geometrically, an element of ${\mathcal F}_n$ is a classical braid  on $n$ strands, with each strand decorated on
the top by an integer, its framing. An element of ${\Bbb Z}^n$, when this is seen as a subgroup of ${\mathcal
F}_n$, is identified with the identity braid on $n$ strands, each strand being decorated by the corresponding
integer of the element.
  For example, the element $f_i$ is  the identity
braid on $n$ strands with framing 1 on the $i$th strand and $0$ elsewhere, while ${f_i}^{a_i}$ is 
 the identity braid on $n$ strands with framing $a_i$ on the $i$th strand and $0$ elsewhere. See Figure~3 for
illustrations.  On the other hand, a braid in $B_n$, when this is seen as a subgroup of ${\mathcal F}_n$, 
is meant as a framed braid with all framings 0. Geometrically, the multiplication in the group ${\mathcal F}_n$ is
the usual concatenation in $B_n$ together with collecting the total framing of each strand to the top of the final
braid. See Figure~4 for an illustration.

\smallbreak 

\begin{figure}[h]
\begin{center}
\includegraphics[width=5.3in]{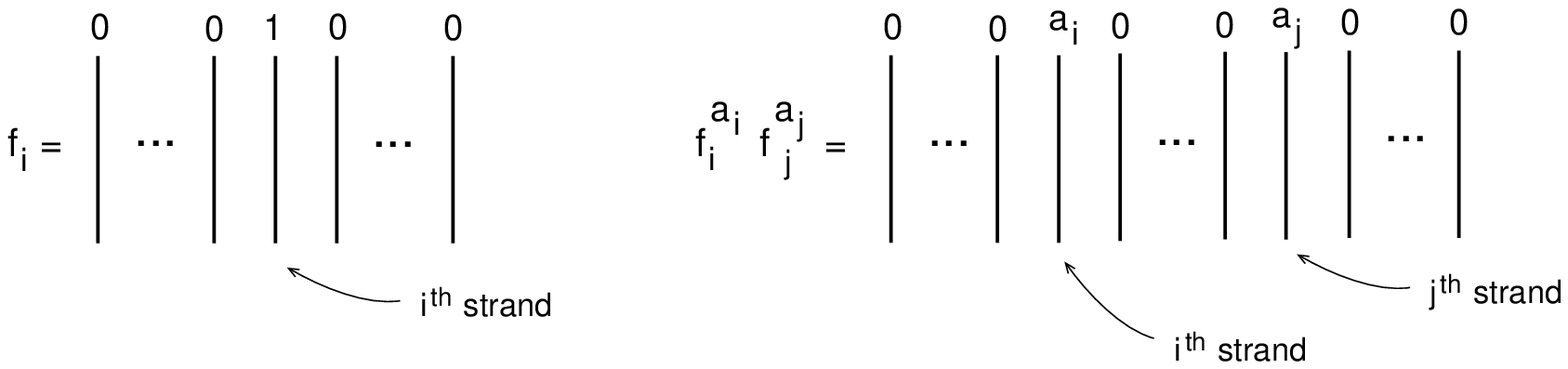}
\end{center}
\caption{Geometric interpretation for $f_i$ and $f_i^{a_i} f_j^{a_j}$}
\label{figure3}
\end{figure}

\smallbreak 

\begin{figure}[h]
\begin{center}
\includegraphics[width=3.7in]{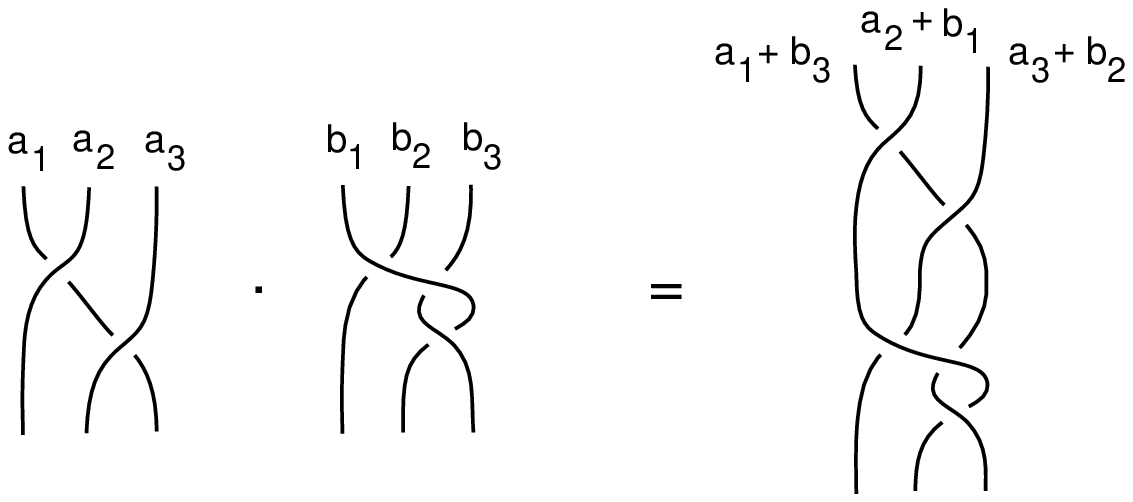}
\end{center}
\caption{Multiplication of framed braids}
\label{figure4}
\end{figure}

\begin{defn}\label{modular} \rm
The {\it  $d$-modular (or simply modular) framed braid group} on $n$ strands 
is defined as   ${\mathcal F}_{d,n} :=  {\Bbb Z}/d{\Bbb Z} \wr B_n = ({\Bbb Z}/d{\Bbb Z})^n \rtimes  B_n$.
\end{defn}

\noindent The group  ${\mathcal F}_{d,n}$ can be considered as the quotient  of ${\mathcal F}_n$ by
imposing the relations
$$
{f_i}^d= 1, \qquad (i =1, \ldots, n).
$$
Clearly, ${\mathcal F}_{d,n}$  has the same geometric interpretation as  ${\mathcal F}_n,$ only that the framings
of the $n$ strands are taken from the cyclic group ${\Bbb Z}/d{\Bbb Z}$. Note now that in ${\mathcal F}_n$ or in
${\mathcal F}_{d,n}$ the $f_i$'s can be deduced from $f_1$, setting for example:
$$
f_i:=\sigma_{i-1} \ldots \sigma_1 f_1 \sigma_1^{-1} \ldots \sigma_{i-1}^{-1}
$$
Then we have the following.

\begin{prop}\label{prefn1} ${\mathcal F}_n$ has a reduced presentation with generators $f_1,
\sigma_1, \ldots, \sigma_{n-1}$ and relations:   
\[
\begin{array}{rcl}
f_1\sigma_j & =  &  \sigma_j f_1 \ \ \ \ \ \ \ \ {\mbox for \  } j>1 \\
f_1\sigma_1 f_1\sigma_1^{-1}  &  =  & \sigma_1 f_1\sigma_1^{-1}f_1    \\  
\sigma_i(\sigma_{i-1}\cdots \sigma_1 f_1\sigma_1^{-1}\cdots 
\sigma_{i-1}^{-1})\sigma_i^{-1}  &  =  & \sigma_i^{-1}(\sigma_{i-1}\cdots 
\sigma_1 f_1 \sigma_{1}^{-1}\cdots \sigma_{i-1}^{-1})\sigma_i \ \ \ {\mbox for \ all \ } i \\  
\end{array}
\]
\noindent and the usual braid relations among the $\sigma_i$'s. 
\end{prop}

\begin{prop}\label{prefn2}  ${\mathcal F}_{d,n}$ has the same
reduced presentation as ${\mathcal F}_n$, but with the extra relation  $f_1^d=1$.
\end{prop}

\subsection{}

In order to define the $p$-adic framed braids we would rather pass to  multiplicative notation for ${\Bbb Z} /
p^r{\Bbb Z}$. Let $C_r$ denote the multiplicative cyclic group of order $p^r$, generated by the element $t_r$. That
is,
$$ 
C_r :=  \langle t_r \  ; \quad t_r^{p^r}  = 1 \rangle.
$$
Then ${\Bbb Z} / p^r{\Bbb Z} \cong C_r$. The maps (\ref{thetarsadd}) of the inverse system  $(C_r , \theta_s^r)$
are now defined by:

\begin{equation}\label{thetars}
\begin{array}{cccc}
\theta_s^r : &  C_r &  \longrightarrow & C_s \\
            \,&  t_r & \mapsto &   t_s
\end{array}
\end{equation}
whenever $r\geq s$. In this notation: $\theta_s^r (t_r^{k_0 + k_1 p + \cdots + k_{r-1} p^{r-1}}) = t_s^{k_0
+ k_1 p + \cdots + k_{s-1} p^{s-1}}$. We have:
$$
{\Bbb Z}_p = \varprojlim C_r
$$ 
and  we can write:
$$
{\Bbb Z}_p =  \{ (t_1^{a_1}, t_2^{a_2}, \ldots) \, \in \prod C_i \ ; \ a_r \in  {\Bbb Z}, \
a_r \equiv a_s \ (\text{mod} \, p^s) \ \ \text{whenever} \ r\geq s\}.
$$ 
For example, for $\underleftarrow{b} = (1, 1+p, 1+p+p^2, \ldots)$ in (\ref{bn}) we write:
\begin{equation}\label{bnmult}
(t_1, t_2^{1+p}, t_3^{1+p+p^2}, \,  \ldots).
\end{equation}
The element 
\begin{equation}\label{fatti}
{\bold t} : =  (t_1, t_2, \ldots) \in \varprojlim C_r 
\end{equation}
corresponds to  $(1, 1, \ldots)$ in the additive notation, so it generates in $\varprojlim C_r$ the constant
sequences.  We shall use the notation ${\Bbb Z} = \langle{\bold t} \rangle$.  
By Corollary~\ref{dense2}, ${\Bbb Z}$ is dense in $\varprojlim C_r$ and ${\bold t}$  is a
topological generator of $\varprojlim C_r$, so an element $(t_1^{a_1}, t_2^{a_2}, \ldots) \in {\Bbb Z}_p$ can be
approximated by the sequence $({\bold t}^{a_i})$ of elements in ${\Bbb Z}$. Following the notation of
(\ref{parrow}) we shall write:
\begin{equation}\label{boldtarrow}
{\bold t}^{\underleftarrow{a}} := (t_1^{a_1}, t_2^{a_2}, \ldots) =  \lim_{k} ({\bold t}^{a_k}).
\end{equation}
For example, the element ${\bold t}^{\underleftarrow{b}}$ in
(\ref{bnmult}) can be approximated by the sequence $({\bold t}, \, {\bold t}^{1+p}, \,  {\bold t}^{1+p+p^2}, \,
\ldots)$:
\begin{equation}\label{limbi}
{\bold t}^{\underleftarrow{b}} = (t_1, t_2^{1+p},t_3^{1+p+p^2}, \ldots) =  \lim_{k} ({\bold
t}^{1+p+\cdots +p^k}).
\end{equation}
With the above notation and according to 1.2, if  ${\bold t}^{\underleftarrow{b}} = (t_1^{b_1}, t_2^{b_2},
\ldots)$ is another element in ${\Bbb Z}_p$, the multiplication in ${\Bbb Z}_p$ is defined as follows:
\begin{equation}\label{pprod}
{\bold t}^{\underleftarrow{a}} {\bold t}^{\underleftarrow{b}} := {\bold t}^{\underleftarrow{a} +
\underleftarrow{b}} = (t_1^{a_1 + b_1}, t_2^{a_2 + b_2}, \ldots)
\end{equation}
and we have the approximation:
\begin{equation}\label{pprodappx}
{\bold t}^{\underleftarrow{a}} {\bold t}^{\underleftarrow{b}} = \lim_{k} ({\bold t}^{a_k+b_k}).
\end{equation}

\subsection{}

Consider now the direct product $C_r^n:= C_r \times\cdots\times C_r$ \, ($n$--times). This is an abelian group 
with the usual product operation defined componentwise, generated by the $n$ elements 
\begin{equation}\label{tiri}
t_{r,i} := (1,\ldots,1, t_r,1,\ldots,1)
\end{equation}
  where $t_r$ is in the $i$th position and where $1$ is the unit
element in $ C_r$. In this notation: 
\begin{equation}\label{tiriprod}
(t_r^{m_1},t_r^{m_2},\ldots,t_r^{m_n}) = t_{r,1}^{m_1} \cdot t_{r,2}^{m_2} \cdots t_{r,n}^{m_n} \ \text{in} \
C_r^n.
\end{equation}
Moreover, $C_r^n$ has the presentation:
\begin{equation}\label{precrn}
C_r^n = \langle t_{r,1}, \ldots,t_{r,n} \  ; \quad t_{r,i}\,t_{r,j} = t_{r,j}\,t_{r,i} \quad \text{and} \quad
t_{r,i}^{p^r}  = (1, \ldots,1) \rangle.
\end{equation}
 Using the maps (\ref{thetars}) of the inverse system  $(C_r , \theta_s^r)$ and (\ref{tiriprod}) we define
componentwise the maps: 
$$
\begin{array}{cccc}
\pi_s^r : &  C_r^n  &  \longrightarrow & C_s^n \\
         &  t_{r,i}  & \mapsto &  t_{s,i}   
\end{array}
$$
whenever $r\geq s$. Then: 
\begin{equation}\label{pirs}
\pi_s^r(t_{r,1}^{m_1} \cdot t_{r,2}^{m_2} \cdots t_{r,n}^{m_n}) = 
t_{s,1}^{m_1(\text{mod} \, p^s)} \cdot t_{s,2}^{m_2(\text{mod} \, p^s)} \cdots t_{s,n}^{m_n(\text{mod} \, p^s)}.
\end{equation}
The maps $\pi_s^r$ are obviously  group epimorphisms, so $(C_r^n , \pi_s^r)$ is an inverse
system of topological groups, indexed by ${\Bbb N}$, and so the inverse limit $\varprojlim C_r^n$ exists. 

\begin{prop}\label{ciarn}
$
\varprojlim C_r^n \cong (\varprojlim C_r)^n = {\Bbb Z}_p^n.
$
\end{prop}
\begin{proof} 
It follows immediately from (\ref{limpro}).
\end{proof} 
Notice now that an element $w\in \varprojlim C_r^n $ can be written as:
\begin{eqnarray*}
w
& = &  
((t_1^{a_{11}},t_1^{a_{12}},\ldots,t_1^{a_{1n}}), \,
(t_2^{a_{21}},t_2^{a_{22}},\ldots,t_2^{a_{2n}}), \ldots )  \\
& = & 
(t_{1,1}^{a_{11}} \, t_{1,2}^{a_{12}} \, \cdots \, t_{1,n}^{a_{1n}}, \, 
t_{2,1}^{a_{21}} \, t_{2,2}^{a_{22}} \, \cdots \, t_{2,n}^{a_{2n}}, \ldots) \quad \text{(by (\ref{tiriprod}))} \\
& = &
 (t_{1,1}^{a_{11}}, t_{2,1}^{a_{21}}, \ldots) \cdot (t_{1,2}^{a_{12}}, t_{2,2}^{a_{22}}, \ldots) \cdots
(t_{1,n}^{a_{1n}}, t_{2,n}^{a_{2n}}, \ldots) \quad \text{(by product operation)}\\
& = &
(t_{r,1}^{a_{r1}} \, t_{r,2}^{a_{r2}} \, \cdots \, t_{r,n}^{a_{rn}})_r .
\end{eqnarray*}
An explicit isomorphism beteween $\varprojlim C_r^n$ and ${\Bbb Z}_p$ is then given by the map:
$$
w\mapsto ((t_1^{a_{11}}, \, t_2^{a_{21}},  \ldots), \, (t_1^{a_{12}}, \, t_2^{a_{22}}, \ldots), 
\ldots, \, (t_1^{a_{1n}}, \, t_2^{a_{2n}}, \ldots))
$$

Thus, we have the identification

\begin{equation}\label{ident1}
(t_{r,1}^{a_{r1}} \, t_{r,2}^{a_{r2}} \, \cdots \, t_{r,n}^{a_{rn}})_r \ 
\stackrel{!}{=}
\  ((t_1^{a_{11}}, \, t_2^{a_{21}},  \ldots), \, (t_1^{a_{12}}, \, t_2^{a_{22}}, \ldots), 
\ldots, \, (t_1^{a_{1n}}, \, t_2^{a_{2n}}, \ldots))
\end{equation}

In particular, the following elements get identified, for $i=1,\ldots,n$: 
$$
\varprojlim C_r^n \ni (t_{r,i})_r \ \stackrel{!}{=} \ ((1, \, 1, \ldots), \ldots, (t_1, \, t_2, \ldots),
\ldots, (1, \, 1, \ldots)) \in (\varprojlim C_r)^n,
$$
where the sequence $(t_1, \, t_2, \ldots)$ is in the $i$th position. Set now 
${\bold 1} := (1,1, \ldots)$ and ${\bold t} = (t_1, \, t_2, \ldots)$ (recall (\ref{fatti})) in $\varprojlim
C_r$ and denote: 
\begin{equation}\label{fattii}
{\bold t}_i := ({\bold 1},\ldots, {\bold 1}, {\bold t} , {\bold 1}, \ldots,{\bold 1} ) \in (\varprojlim C_r)^n
\end{equation}
where  ${\bold t}$  is in the $i$th position. Then we have the
identifications:
\begin{equation}\label{tiri!}
\varprojlim C_r^n \  \ni \ (t_{r,i})_r  \ \stackrel{!}{=} \ {\bold t}_i \ \in \ {\Bbb Z}_p^n
\end{equation}
Thus, with the above notation and with the notation of (\ref{boldtarrow}) we can rewrite the identification
(\ref{ident1}) as follows, for $\underleftarrow{a_i} = (a_{ri})_r$:
\begin{equation}\label{ident2}
\varprojlim C_r^n \ni w = (t_{r,1}^{a_{r1}} \, t_{r,2}^{a_{r2}} \, \cdots \, t_{r,n}^{a_{rn}})_r \ 
\stackrel{!}{=} \ ({\bold t}^{\underleftarrow{a_1}}, {\bold t}^{\underleftarrow{a_2}}, 
\ldots, {\bold t}^{\underleftarrow{a_n}}) = {\bold t}_1^{\underleftarrow{a_1}} {\bold
t}_2^{\underleftarrow{a_2}} \cdots {\bold t}_n^{\underleftarrow{a_n}} \in {\Bbb Z}_p^n.
\end{equation}

\begin{lem}\label{gencrn} 
The identification in $\varprojlim C_r^n$ of the set $X= \{ {\bold t}_1, \ldots, {\bold t}_n \} \subset {\Bbb
Z}_p^n$   is a set of topological generators of \, $\varprojlim C_r^n$. Equivalently, the identification
 in $\varprojlim C_r^n$ of the subgroup ${\Bbb Z}^n = \langle X \rangle$ of ${\Bbb
Z}_p^n$ is dense in $\varprojlim C_r^n$. 
\end{lem}
\begin{proof} By Corollary \ref{dense2} and by Definition \ref{generators},  $ \langle {\bold t}_i\rangle$ is
clearly dense in the $i$th factor 
$(\{ {\bold 1} \}\times \cdots \times \{{\bold 1}\}\times {\Bbb Z}_p \times \{{\bold 1}\}\times \cdots \times
\{{\bold 1}\})$ of  ${\Bbb Z}_p^n$. The result now follows from Corollary \ref{dense2} and the identification
(\ref{ident2}).
\end{proof}
For example, by (\ref{fattii}) and (\ref{tiri!}), and by the approximation (\ref{boldtarrow}), we have the
approximation of $(t_{r,i}^{a_r})_r \in \varprojlim C_r^n$:
\begin{equation}\label{tiriappx}
(t_{r,i}^{a_r})_r  \stackrel{!}{=}  {{\bold t}_i}^{\underleftarrow{a}}  = 
\lim_{k} ({{\bold t}_i}^{a_k}) \stackrel{!}{=} \lim_{k} [(t_{r,i}^{a_k})_r].
\end{equation}
In general, for an element in ${\Bbb Z}_p^n$ we have by (\ref{ident2}), (\ref{boldtarrow}) and (\ref{tiriappx})
the following approximation, where $\underleftarrow{a_i} = (a_{ri})_r$: 
\begin{equation}\label{appxcrn1}
{\Bbb Z}_p^n \ni 
 {\bold t}_1^{\underleftarrow{a_1}}  {\bold t}_2^{\underleftarrow{a_2}} \cdots {\bold
t}_n^{\underleftarrow{a_n}} = \lim_{k} ({\bold t}_1^{a_{k1}}  {\bold t}_2^{a_{k2}} \cdots {\bold
t}_n^{a_{kn}}) = \lim_{k} ({\bold t}^{a_{k1}},{\bold t}^{a_{k2}}, \ldots, {\bold t}^{a_{kn}}).
\end{equation}
Consequently, for the product of two elements in ${\Bbb Z}_p^n$   by (\ref{pprodappx}) the
following approximation, where $\underleftarrow{b_i} = (b_{ri})_r$:
\begin{equation}\label{pnprodappx}
({\bold t}_1^{\underleftarrow{a_1}}  \cdots {\bold t}_n^{\underleftarrow{a_n}}) \, ({\bold
t}_1^{\underleftarrow{b_1}}  \cdots {\bold t}_n^{\underleftarrow{b_n}}) = 
\lim_{k} ({\bold t}_1^{a_{k1}+b_{k1}}  \cdots {\bold t}_n^{a_{kn}+b_{kn}}).
\end{equation}
Hence, for the element in $\varprojlim C_r^n$  $w= (t_{r,1}^{a_{r1}} \,t_{r,2}^{a_{r2}} \, \cdots \,
t_{r,n}^{a_{rn}})_r \stackrel{!}{=} {\bold t}_1^{\underleftarrow{a_1}} {\bold t}_2^{\underleftarrow{a_2}} \cdots
{\bold t}_n^{\underleftarrow{a_n}}$ we obtain, by (\ref{ident2}), (\ref{tiriappx}) and (\ref{appxcrn1}), the
approximation:
\begin{equation}\label{appxcrn2}
\varprojlim C_r^n \ni (t_{r,1}^{a_{r1}} \, t_{r,2}^{a_{r2}} \cdots t_{r,n}^{a_{rn}})_r \  =
 \lim_{k} [(t_{r,1}^{a_{k1}} \, t_{r,2}^{a_{k2}}  \cdots  t_{r,n}^{a_{kn}})_r]
\end{equation}
and for the product of two elements in $\varprojlim C_r^n$ we have the approximation:
\begin{equation}\label{prodappxcrn}
(t_{r,1}^{a_{r1}}  \cdots  t_{r,n}^{a_{rn}})_r \, (t_{r,1}^{b_{r1}} \cdots t_{r,n}^{a_{bn}})_r  =
 \lim_{k} [(t_{r,1}^{a_{k1}+b_{k1}}  \cdots  t_{r,n}^{a_{kn}b_{kn}})_r].
\end{equation}

\subsection{\it $p$-adic framed braids}

In order to introduce the inverse limits in the construction of framed braids we need to start the construction
from the beginning. Consider the cartesian product $C_r^n \times B_n$. Using the maps (\ref{pirs}), define for
any $r, s \in {\Bbb N}$ with
$r\geq s$ the surjective maps: 
\begin{equation}\label{maptimes}
\begin{array}{cccc}
\pi_s^r\times \text{\rm id} : &  C_r^n \times B_n &  \longrightarrow & C_s^n \times B_n \\
         &  (t_{r,1}^{a_{r1}} \, t_{r,2}^{a_{r2}} \, \cdots \, t_{r,n}^{a_{rn}}, \, \sigma)  & \mapsto & 
(t_{s,1}^{a_{s1}} \, t_{s,2}^{a_{s2}} \, \cdots \, t_{s,n}^{a_{sn}}, \, \sigma)   
\end{array}
\end{equation}
 for any $\sigma \in B_n$ and for any exponents satisfying $a_{ri} \equiv a_{si} \
(\text{mod} \, p^s)$. Then we have the following.

\begin{prop}\label{cartesian} 
$(C_r^n \times B_n, \, \pi_s^r\times \text{\rm id})$ is an inverse
system of topological spaces, indexed by ${\Bbb N}$ and we have: 
$$
\varprojlim (C_r^n \times B_n) \cong
\varprojlim C_r^n \times B_n \cong {\Bbb Z}_p^n \times B_n.
$$ 
Moreover, the identification in $\varprojlim (C_r^n \times B_n)$ of ${\Bbb Z}^n \times B_n$ is dense in 
$\varprojlim (C_r^n \times B_n)$ and ${\Bbb Z}^n \times B_n$ is dense in ${\Bbb Z}_p^n \times B_n$.
\end{prop}

\begin{proof} Since the maps $\pi_s^r$ are maps of the inverse system $(C_r^n , \pi_s^r)$, it follows 
immediately that $(C_r^n \times B_n, \, \pi_s^r\cdot \text{\rm id})$ is an inverse system of topological
spaces. An element in $\varprojlim (C_r^n \times B_n)$ is a sequence of the form $((w_1, \sigma),\, (w_2,
\sigma),\, \ldots) $, where  $\sigma \in B_n$ and where $w_1 \in C_1^n, w_2 \in C_2^n, \ldots$, such that
$\pi_s^r(w_r) = w_s$ whenever $r\geq s$. Identifying it with the pair of sequences 
$((w_1, w_2,  \ldots),\, (\sigma, \sigma, \ldots)) \in \varprojlim C_r^n \times \varprojlim B_n$, where
$\varprojlim B_n$ arises as the inverse limit of the trivial inverse system  $(B_n, {\rm id})$, induces the
bijection between $\varprojlim (C_r^n \times B_n)$ and $\varprojlim C_r^n \times  B_n$: 
\begin{equation}\label{ident3}
\varprojlim (C_r^n \times B_n) \ni ((w_1, \sigma),\, (w_2, \sigma),\, \ldots)  \stackrel{!}{=}  
((w_1, w_2,  \ldots), \, \sigma) \in \varprojlim C_r^n \times B_n, 
\end{equation}
where the natural identification between $\varprojlim B_n$  and $B_n$ is induced by the identification $(\sigma,
\sigma,  \ldots) = \sigma$. So the assertion 
$\varprojlim (C_r^n\times B_n) \cong \varprojlim C_r^n \times B_n$ is  proved. Moreover, 
by (\ref{ident1}), $\varprojlim C_r^n \times B_n \cong {\Bbb Z}_p^n \times B_n$.

By  Lemma \ref{gencrn}, and by Corollary \ref{dense2}, the identification of 
${\Bbb Z}^n = \langle {\bold t}_1, {\bold t}_2,\ldots,{\bold t}_n \rangle$ in  $\varprojlim C_r^n$ projects
surjectively on each factor $ C_r^n$ of the inverse system $(C_r^n, \, \pi_s^r)$.  Extending the projection by the
identity map on $B_n$ implies that the identification of ${\Bbb Z}^n\times B_n$ projects surjectively on each
factor $C_r^n \times B_n$ of the inverse system $(C_r^n \times B_n, \pi_s^r\times \text{\rm id})$.  Hence, by 
Corollary \ref{dense1}, the identification of ${\Bbb Z}^n \times B_n$ is dense in  $\varprojlim (C_r^n
\times B_n)$. 
\end{proof}

Consider now the action of the group $B_n$ on the group $C_r^n$ by permutation, as defined in \ref{action}. 
 For the case $d=p^r$ and with the above notation, we have that $C_r^n \rtimes B_n = {\mathcal F}_{p^r,n}$, the
modular framed braid group  with the operation (\ref{frprod}) (in additive notation).
 
\begin{rem} \rm
The generator $f_i$ of
${\mathcal F}_{p^r,n}$ (Proposition \ref{prefn2}) in the additive notation  corresponds to the generator  
$t_{r,i}$ of $C_r^n$. The generators of $C_r^n \rtimes B_n = {\mathcal F}_{p^r,n}$ are the
$n$ elementary framings $t_{r,1}, \ldots, t_{r,n}$ and the $n-1$ elementary braids  $\sigma_1, \ldots,
\sigma_{n-1}$. 
\end{rem}

\noindent Further, use the maps (\ref{maptimes}) of the inverse system  $(C_r^n \times B_n, \, \pi_s^r\times
\text{\rm id})$ to define: 
\begin{equation}\label{mapaction}
\begin{array}{cccc}
\pi_s^r\cdot \text{\rm id} : &  {\mathcal F}_{p^r,n} &  \longrightarrow & {\mathcal F}_{p^s,n} \\
         &  (t_{r,i}, \, id)  & \mapsto &  (t_{s,i}, \, id),   \\
         &((1,\ldots,1), \, \sigma_i) & \mapsto &   ((1,\ldots,1), \, \sigma_i),
\end{array}
\end{equation}
whenever $r\geq s$. 

\begin{lem} $({\mathcal F}_{p^r,n}, \pi_s^r\cdot \text{\rm id})$ is an inverse
system of topological groups, indexed by ${\Bbb N}$. 
\end{lem}

\begin{proof}   On the level of the sets $C_r^n \times B_n$, the map $\pi_s^r\cdot \text{\rm id}$ is
$\pi_s^r\times \text{\rm id}$. We shall show that \ $\pi_s^r\cdot \text{\rm id}$ \ is a group homomorphism.
Indeed, let $(x,\sigma), (y,\tau) \in C_r^n \rtimes B_n$. Then we have:
\begin{eqnarray*}
(\pi_s^r\cdot \text{\rm id}) [(x,\sigma), (y,\tau)] & = & (\pi_s^r\cdot \text{\rm id}) (x \, \sigma(y),\sigma \tau) =
(\pi_s^r(x \, \sigma(y)),\sigma \tau) \stackrel{(in \, C_r^n)}{=} \\
  & = & (\pi_s^r(x) \, \pi_s^r(\sigma(y)),\sigma \tau) \stackrel{(\pi_s^r\circ \sigma=\sigma
\circ\pi_s^r)}{=}  (\pi_s^r(x) \, \sigma(\pi_s^r(y)),\sigma \tau) \\
  & = & (\pi_s^r(x),\sigma)\cdot (\pi_s^r(y) \tau) =  (\pi_s^r\cdot \text{\rm id}) (x,\sigma) \cdot
(\pi_s^r\cdot \text{\rm id}) (y,\tau). 
\end{eqnarray*}

\noindent Hence, $({\mathcal F}_{p^r,n}, \pi_s^r\cdot \text{\rm id})$ is an inverse system of topological groups.
\end{proof}

\begin{defn}\label{finfty}
\rm The {\it $p$-adic framed braid group on $n$ strands} ${\mathcal F}_{\infty ,n}$ is defined to be the inverse
limit of the inverse system $({\mathcal F}_{p^r,n}, \pi_s^r\cdot \text{\rm id})$, that is:
$$
{\mathcal F}_{\infty ,n} := \varprojlim {\mathcal F}_{p^r,n} = \varprojlim (C_r^n \rtimes B_n).
$$
Elements of ${\mathcal F}_{\infty ,n}$ shall be denoted $\underleftarrow{\beta}$.
\end{defn}

\begin{thm}\label{}
There are group isomorphisms: 
$$
{\mathcal F}_{\infty ,n} \cong {\Bbb Z}_p^n \rtimes B_n \cong \varprojlim C_r^n \rtimes B_n. 
$$ 
Moreover, ${\mathcal F}_n$ is dense in ${\Bbb Z}_p^n \times B_n$ and the identification in ${\mathcal F}_{\infty
,n}$ of ${\mathcal F}_n = {\Bbb Z}^n \rtimes B_n$ is dense in ${\mathcal F}_{\infty ,n}$. Finally, 
the identification in ${\mathcal F}_{\infty ,n}$ of the set 
$A = \{ {\bold t}_1,\sigma_1,\ldots, \sigma_{n-1} \} \subset {\mathcal F}_n$
 is a set of topological generators of \, ${\mathcal F}_{\infty ,n}$. 
\end{thm}

\begin{proof}
The second isomorphism is clear from Proposition \ref{ciarn}. We will prove the first one. On the right hand
side $B_n$ acts on ${\Bbb Z}_p^n$ by permutation, that is, a $\sigma \in B_n$ permutes accordingly the
positions of an $n$-tuple of $p$-adic integers. We consider the bijection: 
$$
\alpha :   {\mathcal F}_{\infty ,n}   \longrightarrow  {\Bbb Z}_p^n \rtimes B_n
$$
defined by combining (\ref{ident3}) and (\ref{ident1}). More precisely:
$$
((w_1, \sigma),\, (w_2, \sigma),\, \ldots)   \stackrel{\alpha}{\mapsto}    
([(w_{11}, w_{21}, \ldots), (w_{12}, w_{22}, \ldots),  \ldots, (w_{1n}, w_{2n}, \ldots)], \, \sigma) 
$$
where $w_r = (w_{r1}, w_{r2},\ldots, w_{rn}) \in  C_r^n$.

\smallbreak
\noindent {\it Claim:} $\alpha$ is a group homomorphism. Indeed, let $x = ((w_1, \sigma),\, (w_2, \sigma),\,
\ldots) $ and $y = ((\mu_1, \tau),\, (\mu_2, \tau),\, \ldots) \in {\mathcal F}_{\infty ,n}$, where $\mu_r =
(\mu_{r1}, \mu_{r2},\ldots, \mu_{rn})
\in  C_r^n$. Then:
\begin{eqnarray*}
xy & = & ((w_1, \sigma),\, (w_2, \sigma),\, \ldots) \cdot ((\mu_1, \tau),\, (\mu_2, \tau),\, \ldots) \\
   & = & ((w_1, \sigma) (\mu_1, \tau), (w_2, \sigma) (\mu_2, \tau), \ldots)  \\
   & = & ((w_1 \, \sigma(\mu_1), \sigma\tau), (w_2 \, \sigma(\mu_2), \sigma\tau), \ldots)  \\
   & = & ([(w_{11} \mu_{1\sigma(1)}, \ldots, w_{1n} \mu_{1\sigma(n)}), \sigma\tau], \, 
          [(w_{21} \mu_{2\sigma(1)}, \ldots, w_{2n} \mu_{2\sigma(n)}), \sigma\tau], \ldots). 
\end{eqnarray*}
Hence, 
\begin{eqnarray*}
\alpha(xy) & = & ([(w_{11} \mu_{1\sigma(1)}, w_{21} \mu_{2\sigma(1)}, \ldots), \ldots, 
                   (w_{1n} \mu_{1\sigma(n)}, w_{2n} \mu_{2\sigma(n)}, \ldots)], \sigma\tau).  
\end{eqnarray*}
On the other hand:
\begin{eqnarray*}
\alpha(x) \, \alpha(y) & = & ([(w_{11}, \ldots), \ldots, (w_{1n},  \ldots)], \sigma) \cdot ([(\mu_{11}, 
\ldots), \ldots, (\mu_{1n}, \ldots)], \tau) \\ 
 & = & ([(w_{11},  \ldots), \ldots, (w_{1n},  \ldots)]\, \sigma   
                             ([(\mu_{11}, \ldots), \ldots, (\mu_{1n}, \ldots)], \sigma \tau) \\ 
 & = & ([(w_{11},  \ldots), \ldots, (w_{1n},  \ldots)] \, 
                             [(\mu_{1\sigma(1)}, \ldots), \ldots, (\mu_{1\sigma(n)}, \ldots)], \sigma \tau) \\ 
 & = & ([(w_{11},  \ldots) \, (\mu_{1\sigma(1)}, \ldots), \ldots, (w_{1n},  \ldots) \,  (\mu_{1\sigma(n)},
       \ldots)], \sigma \tau) \\ 
 & = & ([(w_{11} \mu_{1\sigma(1)},\ldots), \ldots, (w_{1n} \mu_{1\sigma(n)}, \ldots) ], \sigma \tau) = \alpha(xy). 
\end{eqnarray*}
Further, ${\Bbb Z}^n \rtimes B_n$ is identical as set  to ${\Bbb Z}^n \times B_n$. By Proposition~\ref{cartesian},
${\Bbb Z}^n \times B_n$ is dense in ${\Bbb Z}_p^n \times B_n$, which in turn is identical as set  to ${\Bbb Z}_p^n
\rtimes B_n$. With similar reasoning the identification in ${\mathcal F}_{\infty ,n}$ of
${\mathcal F}_n = {\Bbb Z}^n \rtimes B_n$ is dense in ${\mathcal F}_{\infty ,n}$.

\smallbreak
For the last statement of Theorem 1, we only need to observe that the generators (\ref{fattii}) of ${\Bbb Z}^n$
are the multiplicative versions of the generators $f_i$ of ${\mathcal F}_n$ given in 2.1. Therefore, the span
$\langle A\rangle$ is isomorphic to the classical framed braid group ${\mathcal F}_n$. So, the identification of
$A$ in ${\mathcal F}_{\infty ,n}$ is a set of topological generators for ${\mathcal F}_{\infty ,n}$.
\end{proof}
\begin{rem}
The fact that ${\Bbb Z}_p$ and $B_n$ contain no elements of finite order imply that ${\mathcal F}_{\infty ,n}
\cong {\Bbb Z}_p^n \rtimes B_n$ contains neither element of finite order. In particular, the modular relations for
the framing are not valid in ${\mathcal F}_{\infty ,n}$. 
\end{rem}

\subsection{\it Geometric interpretations}

By Definition \ref{finfty} a $p$-adic framed braid is an infinite sequence of the same braid $\sigma \in
B_n$, such that the $r$th braid of the sequence  gets framed in the modular framed braid group ${\mathcal
F}_{p^r,n}$ (recall Definition~\ref{modular}) with the framings $(a_{r1}, a_{r2}, \ldots,a_{rn}) \in 
({\Bbb Z}/{\Bbb Z}_{p^r})^n$, where $\underleftarrow{a_i} = (a_{ri})_r$.  In Figure~1 we showed how to interpret a 
 $p$-adic framed braid as an infinite framed cabling of a braid $\sigma \in B_n$. 
 Moreover, by the isomorphism
in Theorem 1, a $p$-adic framed braid can be identified with the element: 
\begin{equation}\label{pfsplit}
{\bold t}_1^{\underleftarrow{a_1}} {\bold t}_2^{\underleftarrow{a_2}} \cdots {\bold t}_n^{\underleftarrow{a_n}}
\cdot \sigma \quad \in {\Bbb Z}_p^n\rtimes B_n,
\end{equation}
that is, the braid $\sigma \in B_n$ with each strand decorated with a $p$-adic integer. See Figure~5.

\smallbreak 

\begin{figure}[h]
\begin{center}
\includegraphics[width=3.7in]{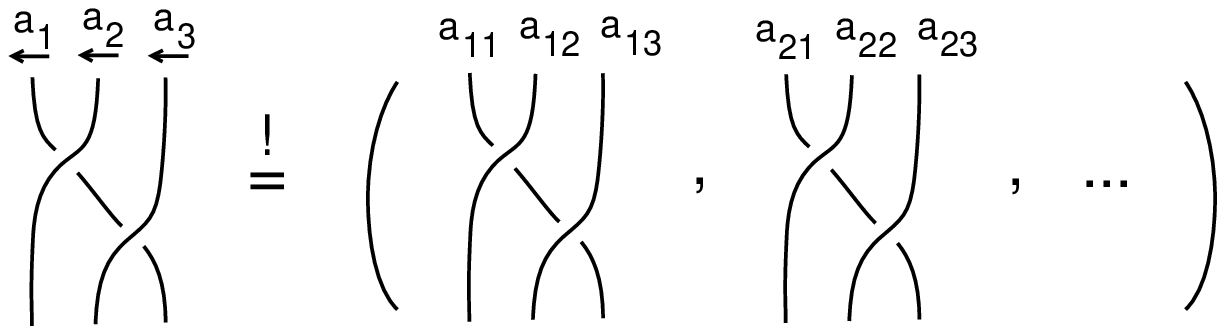}
\end{center}
\caption{A $p$-adic framed braid}
\label{figure5}
\end{figure}

\noindent In particular, the element  ${\bold t}_1^{\underleftarrow{a_1}} {\bold t}_2^{\underleftarrow{a_2}} \cdots
{\bold t}_n^{\underleftarrow{a_n}} \in {\Bbb Z}_p^n$ can be viewed as the identity braid in $B_n$, having the
$p$-adic framing $\underleftarrow{a_i}$ on the $i$th strand, see Figure~6.  

\smallbreak 

\begin{figure}[h]
\begin{center}
\includegraphics[width=1in]{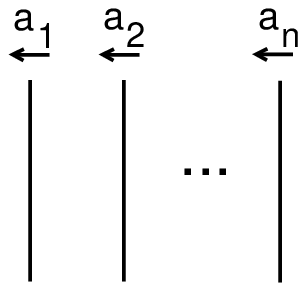}
\end{center}
\caption{A $p$-adic identity framed braid}
\label{figure6}
\end{figure}

\begin{rem}\rm 
As with classical framed braids (\ref{fsplit}), by Theorem 1 and by (\ref{pfsplit}) a $p$-adic framed braid splits
into the 
\lq $p$-adic framing\rq \, part and the \lq braiding\rq \, part. 
\end{rem}

The operation in ${\mathcal F}_{\infty ,n}$ corresponds geometrically to concatenating in each position of the
infinite sequence the two corresponding modular framed braids and collecting the total modular framings to the
top (recall 2.1, (\ref{additive}) and Figure~4). See Figure~7.

\smallbreak 

\begin{figure}[h]
\begin{center}
\includegraphics[width=5in]{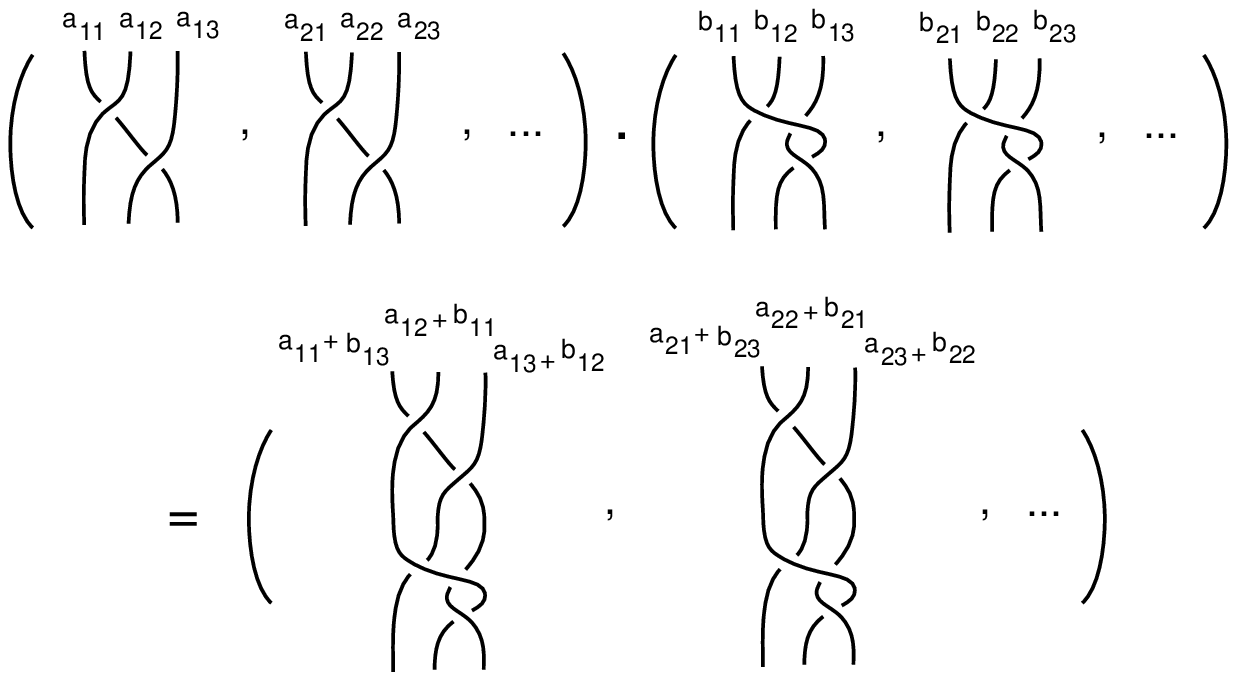}
\end{center}
\caption{Multiplication of $p$-adic framed braids in ${\mathcal F}_{\infty ,n}$}
\label{figure7}
\end{figure}

\noindent On the other hand, by (\ref{pprod}) and 0.3 the multiplication between two elements in  ${\Bbb Z}_p^n
\rtimes B_n$ is defined as follows:
\begin{equation}\label{pbraidprod}
({\bold t}_1^{\underleftarrow{a_1}} \cdots {\bold t}_n^{\underleftarrow{a_n}} \cdot \sigma) \cdot 
({\bold t}_1^{\underleftarrow{b_1}} \cdots {\bold t}_n^{\underleftarrow{b_n}} \cdot \tau) = 
{\bold t}_1^{\underleftarrow{a_1} + \underleftarrow{b_\sigma(1)}} \cdots {\bold t}_n^{\underleftarrow{a_n} +
\underleftarrow{b_\sigma(n)}}
\cdot \sigma\tau
\end{equation}
where  $\underleftarrow{a_i} = (a_{ri})_r$ and $\underleftarrow{b_i} = (b_{ri})_r$. This corresponds  geometrically
to concatenating the two braids $\sigma$ and $\tau$ with $p$-adic framings ($\underleftarrow{a_1},
\ldots, \underleftarrow{a_n}$) and ($\underleftarrow{b_1}, \ldots, \underleftarrow{b_n}$) respectively, and
collecting the total $p$-adic framings to the top. The resulting braid will then have the $p$-adic framings
($\underleftarrow{a_1}+\underleftarrow{b_{\sigma(1)}}, \ldots,
\underleftarrow{a_n}+\underleftarrow{b_{\sigma(n)}}$), where $\underleftarrow{a_i}+\underleftarrow{b_{\sigma(i)}} =
(a_{ri} + b_{r\sigma(i)})_r$, according to (\ref{pprod}). See Figure~8. 

\smallbreak 

\begin{figure}[h]
\begin{center}
\includegraphics[width=3.4in]{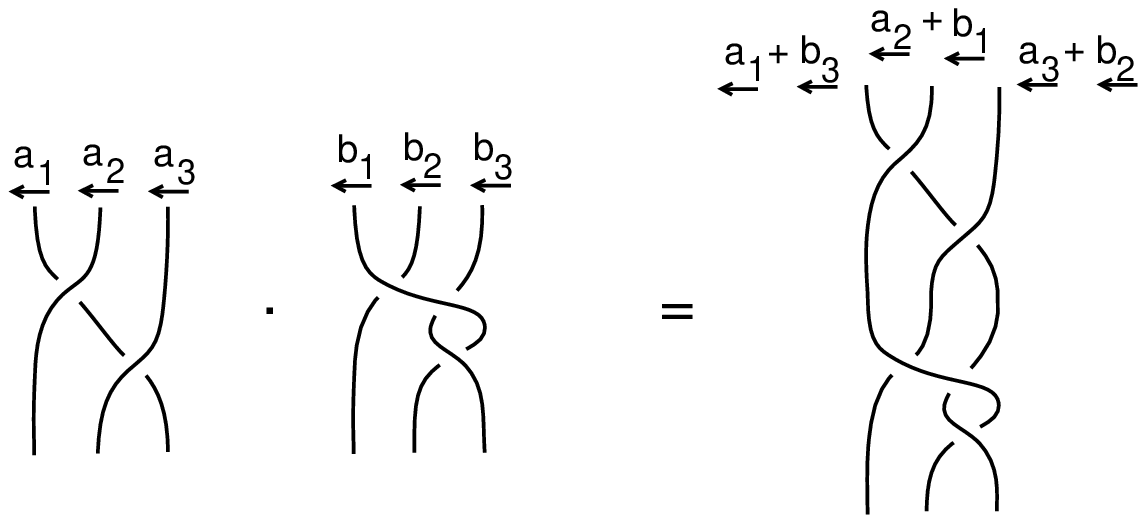}
\end{center}
\caption{Multiplication of $p$-adic framed braids in ${\Bbb Z}_p^n \rtimes B_n$}
\label{figure8}
\end{figure}

In the sequel we will not distinguish between ${\Bbb Z}_p^n \rtimes B_n$ and ${\mathcal F}_{\infty ,n}$, so the
expression (\ref{pfsplit}) and its corresponding geometric interpretation is what we will have in mind from now
on.  In this context, if $\underleftarrow{b} \in {\Bbb Z}_p^n \rtimes B_n$,
such that all framings of $\underleftarrow{b}$ are constant sequences $(k_1),\ldots,(k_n) \in {\Bbb Z}_p^n$ for
$(k_i \in {\Bbb Z})$, then $\underleftarrow{b} \in {\Bbb Z}^n \rtimes B_n$ and it is a classical framed braid with
framings $k_1,\ldots,k_n$. Of course, a classical braid in $B_n$ is meant as a $p$-adic framed braid with all
framings 0.

\subsection{\it Approximations}

By Theorem 1, any element $w= (t_{r,1}^{a_{r1}} \,t_{r,2}^{a_{r2}} \, \cdots \,
t_{r,n}^{a_{rn}} \cdot \sigma)_r$ in ${\mathcal F}_{\infty ,n}$ can be approximated as follows:
\begin{equation}\label{inftyappx}
w  = \lim_{k} (w_k)
\end{equation}
where $w_k$ is the constant sequence $(t_{r,1}^{a_{k1}} \, t_{r,2}^{a_{k2}}  \cdots  t_{r,n}^{a_{kn}} \cdot
\sigma)_r \in {\mathcal F}_{\infty ,n}$. The product of two elements is approximated according to (\ref{inftyappx})
and (\ref{prodappxcrn}).
 
\smallbreak
Further, the fact that ${\mathcal F}_n$ is dense in ${\Bbb Z}_p^n \rtimes B_n \stackrel{!}{=} {\mathcal F}_{\infty
,n}$, means that any $p$-adic framed braid can be approximated by a sequence of classical framed braids. More
precisely, let 
$\underleftarrow{\beta} = {\bold t}_1^{\underleftarrow{a_1}} {\bold t}_2^{\underleftarrow{a_2}}  \cdots {\bold
t}_n^{\underleftarrow{a_n}} \cdot \sigma \in {\Bbb Z}_p^n \rtimes B_n$, where  $\underleftarrow{a_i} = (a_{ri})_r$.
Then,  by (\ref{appxcrn1}), we have:
\begin{equation}\label{pbappx}
\underleftarrow{\beta} =  \lim_k (\beta_k),
\end{equation}
where $\beta_k = {\bold t}_1^{a_{k1}}  {\bold t}_2^{a_{k2}} \cdots {\bold
t}_n^{a_{kn}} \cdot \sigma \, \in {\mathcal F}_n$, and where $a_{ki} \stackrel{!}{=} (a_{ki}, a_{ki},\ldots)$, the
constant sequence in ${\Bbb Z} \subset {\Bbb Z}_p$. 
For example, the $p$-adic braid  ${\bold t}^{\underleftarrow{a}}$, for $\underleftarrow{a} = (a_1,
a_2,\ldots)$, can be approximated as shown in Figure~9, where $a_k \stackrel{!}{=} (a_k, a_k,\ldots) \in {\Bbb
Z} \subset {\Bbb Z}_p$.  See Figure~10 for a generic example. Of course, the product of two
$p$-adic framed braids is approximated accordingly, by (\ref{pbappx}) and (\ref{pnprodappx}).

\smallbreak 

\begin{figure}[h]
\begin{center}
\includegraphics[width=1.5in]{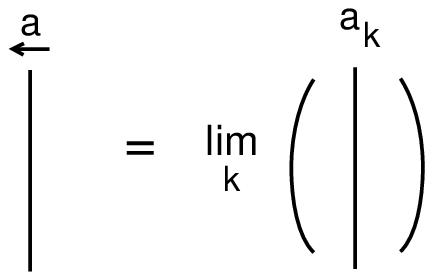}
\end{center}
\caption{The approximation of an one-strand $p$-adic framed braid}
\label{figure9}
\end{figure}

\smallbreak 

\begin{figure}[h]
\begin{center}
\includegraphics[width=4.1in]{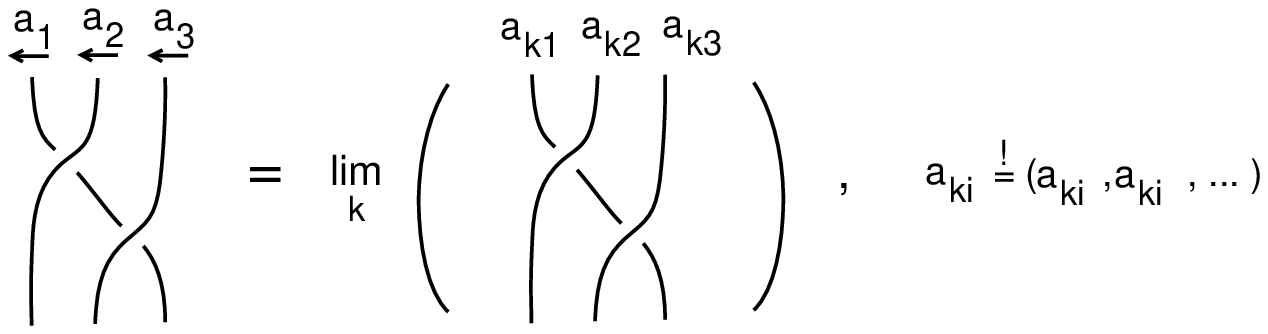}
\end{center}
\caption{The approximation of a $p$-adic framed braid}
\label{figure10}
\end{figure}

\section{Quotient algebras from $p$-adic framed braids}

In this section we define the main algebra studied in the paper. This algebra arises as the 
inverse limit of an inverse system of so-called Yokonuma-Hecke algebras. In the
sequel we fix an  element $u$  in ${\Bbb C}\backslash \{0\}$ and we shall denote 
${\Bbb C}[G]$ (or simply ${\Bbb C}G$) the group algebra of a group $G$. 

\subsection{} 
 
Let $H=\langle h \rangle$ be a  finite cyclic group of order $d$. As in (\ref{tiri})  we define
the element $h_i$ in $H^n: = H\times\cdots \times H$ ($n$ copies) as the
element having $h$ on the $i$th component and $1$ elsewhere. So, for any element 
 $(h^{a_1},\ldots , h^{a_n})\in H^n$ we have
$$
(h^{a_1},\ldots , h^{a_n}) = h_1^{a_1}\cdots h_n^{a_n}.
$$
 For any $i, j$ with $i\not=j$, we define the subgroups $H_{i,j}$  of $H^n$ as follows:
\begin{equation}\label{hij}
H_{i, j} := \langle h_ih_j^{-1} \rangle.
\end{equation}
Clearly, $H_{i,j}$ is isomorphic to  the group $H$. In ${\Bbb C}[H^n] = {\Bbb C}H^n$ we define
the following  elements:
$$
e_{d,i,j} := \frac{1}{d} \sum_{x\in  H_{i,j}}x \in {\Bbb C}H^n
$$
or, equivalently:
$$
e_{d,i,j} = \frac{1}{d} \sum_{1\leq m \leq d}h_i^mh_{j}^{-m}.
$$
\begin{lem}\label{ideeri}
For any $i, j$ with $i\not= j$ the elements $e_{d,i,j}$ are idempotents.
\end{lem}

\begin{proof}
It is enough to observe that $e_{d,i,j}$ is the average on the 
elements of the group $H_{i,j}$. Indeed,
$$
(e_{d,i,j})^2 =  \frac{1}{d} \sum_{y\in  H_{i,j}}y \frac{1}{d} \sum_{x\in  H_{i,j}}x 
=  \frac{1}{d^2} \sum_{y\in  H_{i,j}} \sum_{x\in  H_{i,j}}yx = 
 \frac{d}{d^2} \sum_{x^{\prime}\in  H_{i,j}}x^{\prime} = e_{d,i,j}.
$$
\end{proof} 

\begin{rem} \rm
Notice that $H_{i,j}=H_{j,i}$. In the case $j=i+1$ we
denote  $H_{i,i+1}$ by $H_{i}$ and $e_{d, i, i+1}$ by $e_{d, i}$.
\end{rem}

\subsection{}

Consider now the modular framed braid group ${\mathcal F}_{d,n}$ (Definition \ref{}). 
The ${\Bbb C}$--algebra ${\Bbb C}H^n$ is a subalgebra of the group algebra ${\Bbb C}{\mathcal
F}_{d,n}$ and 
 the elements $e_{d,i,j}$ are still idempotents in ${\Bbb C}{\mathcal F}_{d,n}$. 
 The main commutation relations among them and the elementary braids $\sigma_i$ are given
   in the proposition below.
\begin{prop} \label{peri}
For any $i,j \in \{1, \ldots , n-1\}$ we have:
\begin{enumerate}
\item $\sigma_i^{\pm 1}e_{d, j}=
e_{d,j}\sigma_i^{\pm 1}$, for all $j\not= i-1, i+1$
\item $\sigma_i^{\pm 1} e_{d ,j} =
e_{d , i,j}\sigma_i^{\pm 1}$,  for $\vert i-j\vert =1$
\item $e_{d,j} \sigma_i^{\pm 1}  =
\sigma_i^{\pm 1}e_{d , i, j}$,  for $\vert i-j\vert =1$
\item $e_{d,i}\, h_1^{a_1}\cdots h_n^{a_n} = e_{d,i}\, h_1^{a_1}\cdots h_{i-1}^{a_{i-1}}
 \left( h_i^{a_{i+1}}h_{i+1}^{a_i}\right) h_{i+2}^{a_{i+2}}\cdots h_n^{a_n}.$
\end{enumerate}
\end{prop}

\begin{proof}
 (1) If $j\not = i, i\pm 1$, the claim follows from the fact that $\sigma_i$ commutes with  $h_j$. 
 Let now that $j= i$. We have  $\sigma_ie_{d,i} = \sigma_i d^{-1} \sum_{s} h_i^s h_{i+1}^{-s}.$ Note
now that
 $\sigma_i h_i^s h_{i+1}^{-s}=h_{i+1}^s\sigma_ih_{i+1}^{-s} =
 h_{i+1}^sh_i^{-s}\sigma_i.$ Then
$$
\sigma_ie_{d,i} = \frac{1}{d} (\sum_{s}h_{i+1}^sh_i^{-s})\sigma_i=e_{d,i}\sigma_i.
$$

\noindent (2) Let $j = i +1$. We have that  
$\sigma_ih_{i+1}^sh_{i+2}^{-s} = h_i^s\sigma_ih_{i+2}^{-s} =
h_i^s f_{i+2}^{-s}\sigma_i.$ So, we deduce: $\sigma_ie_{d,i+1} =
d^{-1} \sum_{s} h_i^s h_{i+2}^{-s}\sigma_i.$
 Claim (3) follows in the same way as (2).

\smallbreak
\noindent (4)  Setting $c:= h_1^{a_1}\cdots h_n^{a_n}$  we have:
\begin{eqnarray*}
h_i^s h_{i+1}^{-s}c 
& = &
h_1^{a_1}\cdots h_{i-1}^{a_{i-1}}h_i^{a_i+s}h_{i+1}^{a_{i+1}-s}h_{i+2}^{a_{i+2}}\cdots
h_n^{a_n}\\
 & = &
h_1^{a_1}\cdots h_{i-1}^{a_{i-1}}h_i^{(s+
a_i-a_{i+1})+a_{i+1}}h_{i+1}^{-(s+ a_i -a_{i+1})+a_i}h_{i+2}^{a_{i+2}}\cdots h_n^{a_n}\\
 & = &
h_1^{a_1}\cdots h_{i-1}^{a_{i-1}}h_i^{(s+
a_i-a_{i+1})}h_i^{a_{i+1}}h_{i+1}^{-(s+ a_i -a_{i+1})}h_{i+1}^{a_i}h_{i+2}^{a_{i+2}}\cdots
 h_n^{a_n}\\
& = &
(h_i^{(s+a_i-a_{i+1})}h_{i+1}^{-(s+ a_i -a_{i+1})})
h_1^{a_1}\cdots h_{i-1}^{a_{i-1}}h_i^{a_{i+1}}h_{i+1}^{a_i}h_{i+2}^{a_{i+2}}\cdots h_n^{a_n}.
\end{eqnarray*}

Therefore,

\begin{eqnarray*}
e_{d,i} c &  = & \frac{1}{d} \sum_{0\leq s\leq d-1} h_i^s h_{i+1}^{-s} c \\
&  = &
\left( \frac{1}{d}  \sum_{s}h_i^{(s+a_i-a_{i+1})}h_{i+1}^{-(s+ a_i -a_{i+1})}\right) \\
& &
h_1^{a_1}\cdots h_{i-1}^{a_{i-1}}\left(h_i^{a_{i+1}}h_{i+1}^{a_i}\right)h_{i+2}^{a_{i+2}}\cdots
h_n^{a_n} \\
& = &
e_{d,i} h_1^{a_1}\cdots
h_{i-1}^{a_{i-1}}\left(h_i^{a_{i+1}}h_{i+1}^{a_i}\right) h_{i+2}^{a_{i+2}}\cdots h_n^{a_n}.
\end{eqnarray*}
\end{proof}
\begin{rem}\rm
The elements $h_i$ correspond to the elementary framings $f_i$ in the additive notation of 
subsection 2.1 and, for $d=p^r$, to the elements $t_{r,i}$ defined in (\ref{tiri}).
\end{rem}

\subsection{}

The Yokonuma--Hecke (Y--H) algebras were introduced by T.~Yokonuma \cite{yo}
in the context of  Chevalley groups, as   generalizations 
 of the Iwahori-Hecke algebras. More precisely, the
 Iwahori-Hecke algebra associated to a finite Chevalley group $G$ is the
 centralizer algebra associated to the permutation representation of $G$ 
 with respect to a 
 {\it Borel subgroup} of $G$. The Y--H algebra is the centralizer 
  algebra associated to the permutation representation of $G$ with respect 
  to  a {\it unipotent subgroup} of $G$. So, the  Y--H algebra can be 
 also regarded  as a particular  case of a unipotent algebra. See \cite{th} for the general  
definition of unipotent algebras.

\smallbreak

We define the {\it Yokonuma--Hecke algebra of type} $A$, ${\rm Y}_{d, n}(u)$, as the quotient of the group
algebra of the modular framed braid group ${\mathcal F}_{d,n}$ under the quadratic relations: 
\begin{equation}\label{ceri}
g_i^2 =  1 + (u-1) \, e_{d,i} \, (1-g_i) \qquad (i=1,\ldots , n-1).
\end{equation}
More precisely, 
 ${\rm Y}_{d, n}(u)$ is defined as follows:
$$
{\rm Y}_{d, n}(u) : =\frac{{\Bbb C}{\mathcal F}_{d,n}}
 {\langle \sigma_i^2 - 1 - (u-1)e_{d,i}(1-\sigma_i),
 \quad i=1,\ldots , n-1 \rangle}.
$$

Corresponding now $\sigma_i \in {\Bbb C}{\mathcal F}_{d,n}$ to $g_i \in {\rm Y}_{d,
n}(u)$ and $f_i\in {\mathcal F}_{d,n}$ to $h_i\in {\rm Y}_{d, n}(u)$, we obtain from 
the above and from Proposition~\ref{prefn2} a reduced presentation of ${\rm Y}_{d, n}(u)$, by setting:
\begin{equation}\label{tiis}
h_i= g_i \cdots g_1 h_1 g_1^{-1} \cdots g_i^{-1}.
\end{equation}
Then a reduced presentation of ${\rm Y}_{d, n}(u)$ is given in the theorem below.
\begin{thm}\label{yrnr}
The algebra ${\rm Y}_{d,n}(u)$ can be presented with the generators 
 $h_1$, $g_1$, $\ldots,$
 $g_{n-1}$ and the following relations:
\begin{enumerate}
\item Braid relations among the $g_i$'s
\item $h_1 g_i =g_i h_1$, for $i\geq 2$
\item $h_1g_1h_1 g_1^{-1} =  g_1h_1 g_1^{-1}h_1$
\item $h_1^{d}=1$
\item
$g_i(g_{i-1}\cdots g_1h_1g_1^{-1}\cdots g_{i-1}^{-1})g_i^{-1}
= g_i^{-1}(g_{i-1}\cdots g_1h_1g_1^{-1}\cdots
g_{i-1}^{-1})g_i$
\item $g_i^2 = 1 + (u-1)e_{d,i} (1-g_i)$,

$(i=1,\ldots , n-1).$
\end{enumerate}
\end{thm}

In this above notation, we may rewrite the elements $e_{d,i}\in {\rm Y}_{d, n}(u)$ as:
$$
e_{d,i} = \frac{1}{d}
\sum_{1\leq m\leq d}(g_{i-1}^{-1}\cdots g_1^{-1}h_1^mg_1\cdots
g_{i-1}) (g_i\cdots g_1h_1^{-m}g_1^{-1}\cdots g_i^{-1}).
$$
\begin{rem} \rm
The  Y--H algebra ${\rm Y}_{d, n}(u)$ can be also thought of as a $u$--deformation of the group algebra
${\Bbb C}[H\wr S_n]$ in the following sense: The  algebra ${\Bbb C}[H\wr S_n] =
{\Bbb C} [ H^n \rtimes S_n]$ contains  ${\Bbb C}H^n$ as a subalgebra, so the elements 
$e_{d,i}$ are also in ${\Bbb C}[H\wr S_n]$. We correspond now the generator
$s_i \in {\Bbb C}[H\wr S_n]$ to the generator $g_i \in {\rm Y}_{d, n}(u)$, the  generator $h_1 \in
{\Bbb C}[H\wr S_n]$ to the  generator $h_1\in {\rm Y}_{d, n}(u)$ and $e_{d,i} \in {\Bbb C}[H\wr S_n]$ to
$e_{d,i} \in {\rm Y}_{d, n}(u)$ (we keep the same notation).
 Then, the canonical presentation of ${\Bbb C} [H\wr S_n]$ gives rise to a
  presentation  of ${\rm Y}_{d, n}(u)$ (the same as in Theorem \ref{yrnr}) by imposing the quadratic
relations in (\ref{ceri}) instead of the relations $s_i^2=1$.
\end{rem}

\begin{rem}\rm
The fact that the element $e_{d,i}$ is an idempotent makes it possible to define in ${\rm Y}_{d,n}(u)$
the inverse of $g_i$. Indeed,  multiplying relation (\ref{ceri}) by $g_i$ gives  $g_i^3 =
g_i + (u-1) \, e_{d,i}g_i - (u-1) \, e_{d,i} \, g_i^2$. Replacing now $g_i^2$ by its
expression (\ref{ceri}) and using the fact that $e_{d,i}$ is an idempotent, we obtain that $g_i^3
= g_i -(u^2 - u )e_{d,i} + (u^2 - u)e_{d,i}g_i$. Using again (\ref{ceri}) we  substitute  
$e_{d,i}g_i$ by $(u-1)^{-1}(1 + (u-1)e_{d,i} - g_i^2)$, so we have $ g_i^3 = u + g_i - 
u g_i^2$. Multiplying the latter by $g_i^{-1}$ we deduce  $g_i^{-1} = u^{-1}( g_i^2 + ug_i - 1)$ and, 
 using again (\ref{ceri}), we finally obtain:
\begin{equation}\label{invrs}
g_i^{-1} = g_i - (u^{-1} - 1)\, e_{d,i} + (u^{-1} - 1)\, e_{d,i}\, g_i.
\end{equation}
\end{rem}

\subsection{}

In this part we give  a   diagrammatic interpretation of the elements $e_{d,i}$  and of the quadratic
relations in ${\rm Y}_{d, n}(u)$. The elements $e_{d,i}$ seen  as elements of ${\Bbb C}{\mathcal
F}_{d,n}$   can be interpreted geometrically as the average of the sum of $d$ identity framed braids with 
framings as shown in Figure~11.

\bigbreak 
\begin{figure}[h]
\begin{center}
\includegraphics[width=3.2in]{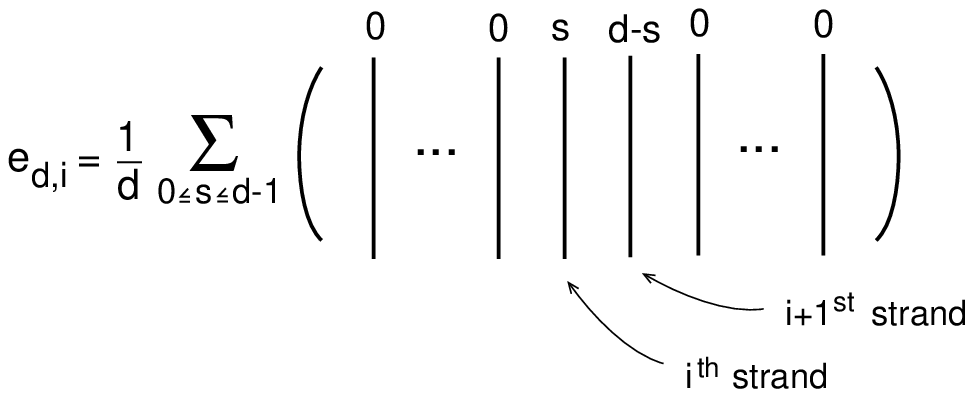}
\end{center}
\caption{The elements $e_{d,i}$}
\label{figure11}
\end{figure}

Similarly, the quadratic relations $g_i^2 =  1 + (u-1) \, e_{d,i}\cdot 1 - (u-1)e_{d,i}g_i$ can be also
considered as relations in ${\Bbb C}{\mathcal F}_{d,n}$. In Figure~12 we illustrate the relation for
$g_1^2$ in  ${\Bbb C}{\mathcal F}_{d,3}$. Note that the effect of
$e_{d,i}$ on the identity element or on  $g_i$ is to produce $d$ copies and frame appropriately the $i$th
and $(i+1)$st strand. Similar is the effect of $e_{d,i}$ on any braid. In Figure~13 we illustrate the
quadratic relation in a compact form.  Finally, in Figure~14 we illustrate the equation for
$\sigma_1^{-1}$ in ${\Bbb C} {\mathcal F}_{d,3}$.

\bigbreak 
\begin{figure}[h]
\begin{center}
\includegraphics[width=5.2in]{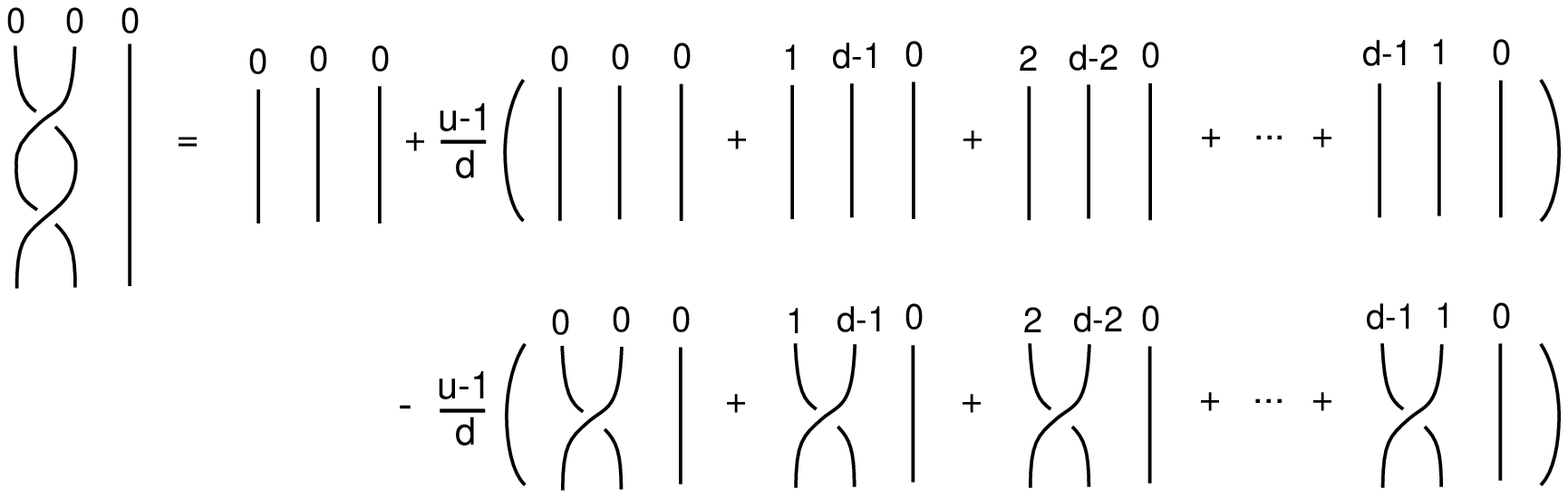}
\end{center}
\caption{Geometric interpretation of $g_1^2$}
\label{figure12}
\end{figure}

\bigbreak 
\begin{figure}[h]
\begin{center}
\includegraphics[width=3.5in]{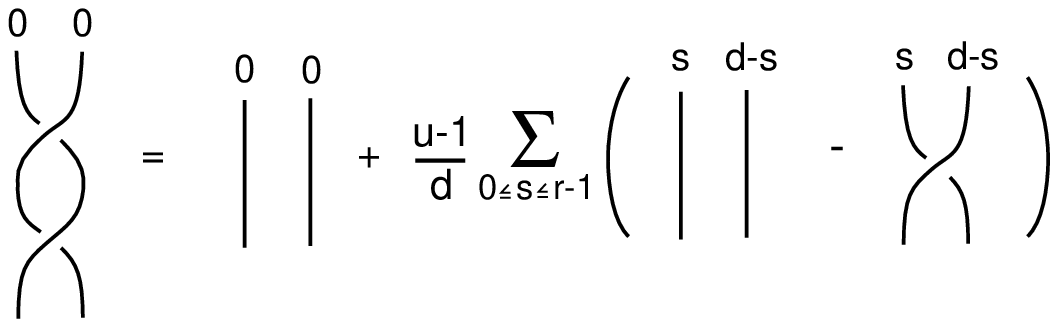}
\end{center}
\caption{$g_i^2 = 1 + (u-1) \, e_{d,i} \, (1-g_i)$}
\label{figure13}
\end{figure}

\bigbreak 
\begin{figure}[h]
\begin{center}
\includegraphics[width=5.2in]{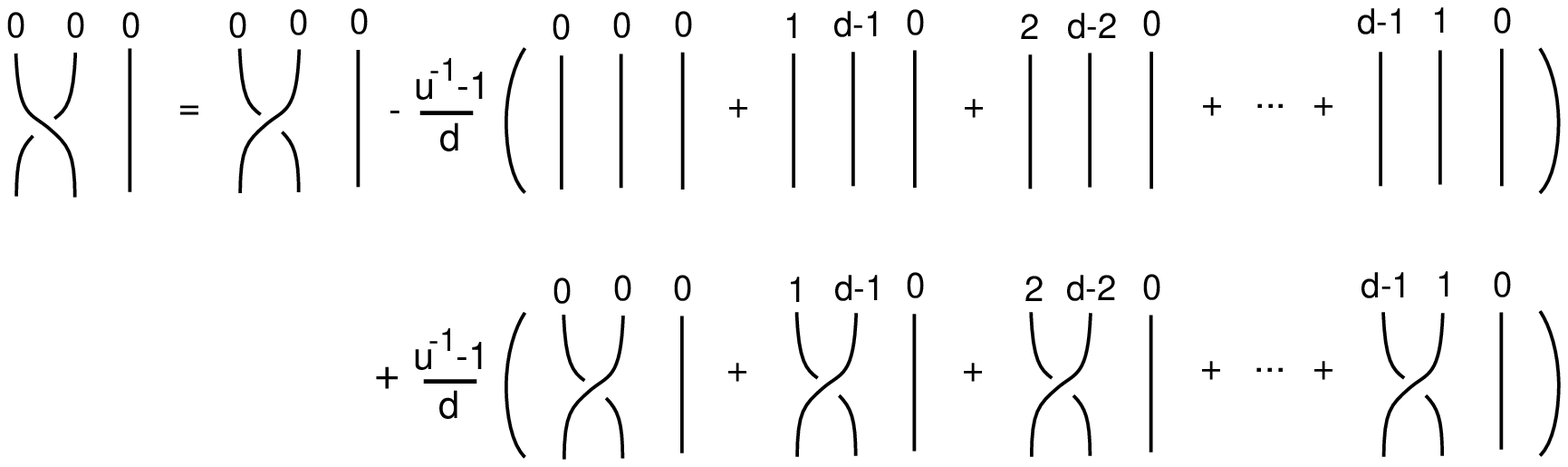}
\end{center}
\caption{Geometric interpretation of $g_1^{-1}$}
\label{figure14}
\end{figure}

\begin{rem}\rm
Note the resemblance of this relation to the skein relations used for defining classical
quantum link invariants.  For  $d=1$ the relation becomes the skein relation of the 2-variable
Jones polynomial (HOMFLYPT), that arises from the quadratic relation of the Hecke algebra of type $A$, 
see \cite{jo}. In fact, ${\rm Y}_{1, n}(u)$ coincides with the Hecke algebra of type $A$.
\end{rem}

\subsection{\it The $p$-adic Yokonuma--Hecke algebra}

 We shall now explain our construction of the $p$-adic Yokonuma--Hecke algebra   
 ${\rm Y}_{\infty, n}(u)$. The ${\Bbb C}$--algebra  ${\rm Y}_{\infty, n}(u)$  will be defined as
the inverse limit of an inverse system of the  Y--H algebras ${\rm Y}_{p^r, n}(u), \, r\in {\Bbb N}$,
wherer $p$ is a fixed prime number. On this family of Y--H algebras we  consider  epimorphisms 
$$
\varphi_s^r: {\rm Y}_{p^r, n}(u)\longrightarrow {\rm Y}_{p^s, n}(u)\qquad (r\geq s),
$$ 
induced from the group homomorphisms $\pi_s^r\cdot {\rm id}$ defined in (\ref{pirs}). More
precisely, extending $\pi_s^r\cdot {\rm id}$ linearly, yields a natural algebra epimorphism
$$
\phi_s^r: {\Bbb C}{\mathcal F}_{p^r,n}\longrightarrow {\Bbb C}{\mathcal F}_{p^s, n}\qquad (r\geq
s).
$$ 

It is a routine to check the following lemma.
\begin{lem}
 $({\Bbb C} {\mathcal F}_{p^r,n}, \phi_s^r)$ is an  inverse system of rings, indexed 
by ${\mathbb N}$.
\end{lem}

Note that the natural embedding $\iota_r: {\mathcal F}_{p^r,n} \hookrightarrow {\Bbb
C}{\mathcal F}_{p^r,n}$  induces a natural embedding 
$\varprojlim \iota_r: {\mathcal F}_{\infty,n} \hookrightarrow \varprojlim {\Bbb C}{\mathcal
F}_{p^r,n}$. So, up to  identifications, we have the inclusions:
$$
{\mathcal F}_n \subseteq {\mathcal F}_{\infty,n}\subseteq   \varprojlim {\Bbb C} {\mathcal
F}_{p^r,n}
$$
Recall now that ${\bold t}_1 :=({\bold t}, {\bold 1}, \ldots , {\bold 1})$ and 
 $\sigma_i :=(\sigma_i, \sigma_i, \ldots )$ in $\varprojlim {\Bbb C}{\mathcal
F}_{p^r,n}$, we have the following result:

\begin{prop}\label{glga}
The set  $X=\{{\bold t}_1, {\sigma}_1, \ldots, {\sigma}_{n-1}\}$ is a
 set  of topological  generators of the algebra  $\varprojlim {\Bbb C}{\mathcal
F}_{p^r,n}$. In particular, the subalgebra ${\Bbb C}{\mathcal F}_n$ is dense  in
$\varprojlim {\Bbb C}{\mathcal F}_{p^r,n}.$
\end{prop}

\begin{proof} 
By Proposition \ref{prefn1}, the set $X$ is
a set of generators for  the group ${\mathcal F}_n$, hence  $X$ spans the
 algebra ${\Bbb C}{\mathcal F}_n$.   Now, the mapping $\sigma_i \mapsto \sigma_i$,
${\bold t}_1\mapsto t_{r,1}$ defines an epimorphism
$\eta_r: {\Bbb C}{\mathcal F}_n \longrightarrow  {\Bbb C}{\mathcal F}_{p^r,n}$,
for  any $r\in {\Bbb N}$. Notice now that $\eta_r$ is surjective and that
we have the following commutative diagram:


$$
\begin{CD}
1@>>> {\Bbb C}{\mathcal F}_n  @>>>  \varprojlim {\Bbb C}{\mathcal F}_{p^r,n} \\
 & \, &  \eta_r @VVV      \xi _r @VVV       \\
 & \, & {\Bbb C}{\mathcal F}_{p^r, n}  @= {\Bbb C}{\mathcal F}_{p^r, n}
\end{CD}
$$
where  $\xi_r$ is the natural projection. Then the proof follows from Corollary~\ref{dense1}.
\end{proof}

Recall now the subgroups $H_{i,j}$ defined in (\ref{hij}). With the notations 
of Section~2 for $H=C_r$ we denote these subgroups  by $H_{r,i,j}$ and we have:  
$$
H_{r,i,j} = \langle t_{r,i}t_{r,i+1}^{-1}\rangle.
$$
Hence 
$e_{p^r,i,j}\in  {\Bbb C}C_r^n$. Recalling also that ${\mathcal F}_{p^r,n} = 
C_r^n \rtimes B_n$, we have the following.

\begin{prop}\label{proeri}
For any $i, j$ with $i\not= j$ and for $s\leq r$, we have:
\begin{enumerate}
\item The homomorphism $\phi_s^r$ maps $H_{r,i,j}$ onto $H_{s,i,j}$.
\item The  kernel of the restriction of $\phi_s^r$ on $H_{r,i,j}$ has order $p^{r-s}$.
\item $\phi_s^r(e_{p^r,i,j}) = e_{p^s,i,j}$.
\end{enumerate}
\end{prop}
\begin{proof}
Since $\phi_s^r(t_{r,i}t_{r,j}^{-1}) = t_{s,i}t_{s,j}^{-1}$ claim (1) follows. Claim (2) is clear by the
fundamental theorem of homomorphisms for groups. Finally, claim (3) follows directly from claims (1) and
(2).
\end{proof}

Defining now in ${\Bbb C}{\mathcal F}_{p^r,n}$ the elements: 
$$
\epsilon_{r,i} := \sigma_i^2 - 1 - (u-1)e_{p^r,i}(1-\sigma_i) \in 
{\Bbb C}{\mathcal F}_{p^r,n}\qquad (i=1,\ldots , n-1) ,
$$
and the ideal 
$$
I_{p^r, n} = \langle \epsilon_{r,i} \, ;\,i=1,\ldots , n-1 \rangle.
$$
 We have that 
 $$
 {\rm Y}_{p^r, n}(u) =\frac{{\Bbb C}{\mathcal F}_{p^r,n}}{I_{p^r, n}}.
 $$
Using (3) of Proposition \ref{proeri} we obtain the following lemma.

\begin{lem}\label{ideiri}
For all $i$ and for $s\leq r$, we have:  $\phi_s^r(I_{p^r,n})=  I_{p^s,n}.$
\end{lem}

According to Lemma \ref{ideiri} we obtain the following commutative diagram of
rings:
$$
\begin{CD}
{\Bbb C}{\mathcal F}_{p^r, n}  @>\phi_s^r>>{\Bbb C} {\mathcal F}_{p^s,n} \\
   @V \rho_rVV         @V\rho_s VV    \\
{\rm Y}_{p^r, n}(u) @> \varphi_s^r  >>  {\rm Y}_{p^s,n}(u)\\
\mbox{}
\end{CD}
$$
where $\rho_r$ and $\rho_s$ are the canonical epimorphisms and $\varphi_s^r$ 
is defined  via $\phi_s^r$ as: 
$$
\varphi_s^r (x + I_{p^r , n }) : = \phi_s^r (x) + I_{p^s , n }.
$$
Recall that ${\rm Ker}(\rho_r) = I_{p^r,n}.$ Thus,  the inverse system
 $\left( {\Bbb C}{\mathcal F}_{p^r,n}, \phi_s^r\right)$ induces  the inverse system
$$
\left({\rm Y}_{p^r,n}(u), \varphi_s^r\right),
$$ 
indexed by ${\Bbb N}$.

\begin{defn}\label{yinfty}
The  {\it $p$-adic Yokonuma--Hecke algebra } ${\rm Y}_{\infty ,n}(u)$ is defined as the inverse limit of
this last inverse system.
$$
{\rm Y}_{\infty ,n}(u): = \varprojlim {\rm Y}_{p^r,n}(u).
$$
\end{defn}
The algebra   ${\rm Y}_{\infty ,n}(u)$ is equipped  with canonical epimorphisms:
$$
\Xi_r: {\rm Y}_{\infty ,n}(u)\longrightarrow {\rm Y}_{p^r ,n}(u),
$$
 such that $\varphi_s^r\circ \Xi_r = \Xi_s$.

\subsection{}

We shall now try to understand better the structure of ${\rm Y}_{\infty ,n}(u)$. 
By Lemma \ref{ideiri} the
restriction of $\phi_s^r$ to $I_{p^r, n}$ yields the  inverse system 
 $\left( I_{p^r,n}, \phi_s^r\right)$.
Furthermore, for any $r$ we have the following exact sequence:
$$
\begin{CD}
0 @>>>I_{p^r, n}  @>\iota_r>> {\Bbb C}{\mathcal F}_{p^r,n} @>\rho_r>>{\rm Y}_{p^r ,n}(u)@>>> 0
\end{CD}
$$
Then, by  (\ref{exact}), we obtain the exact sequence:
$$
\begin{CD}
0 @>>>\varprojlim I_{p^r, n}  @>\iota>>  \varprojlim{\Bbb C}{\mathcal F}_{p^r,n}
 @>\rho>>  {\rm Y}_{\infty ,n}(u)
\end{CD}
$$
where $\iota:= \varprojlim \iota_r$ and $\rho:= \varprojlim\rho_r$. Hence, 
 and since $\varprojlim I_{p^r, n}$ is an ideal in $\varprojlim{\Bbb C}{\mathcal F}_{p^r,n}$, 
we have:
$$
\frac{\varprojlim{\Bbb C}{\mathcal F}_{p^r,n}}{\varprojlim I_{p^r, n}} \cong
\rho ( \varprojlim{\Bbb C}{\mathcal F}_{p^r,n}).
$$
At this writing it is not clear whether the map $\rho$ is a surjection or not. Yet, we have the following
result.

\begin{prop}\label{prodense}
$\rho ( \varprojlim{\Bbb C}{\mathcal F}_{p^r,n})$ is dense in
 $ {\rm Y}_{\infty ,n}(u)$.
\end{prop}
\begin{proof}
The proof is again an application of Corollary \ref{dense1}. Indeed, define
the map $\theta: \rho (x) \mapsto (\rho_r \circ \xi_r)(x)$, for
$x = (x_r)\in \varprojlim{\Bbb C}{\mathcal F}_{p^r,n}$. Clearly $\theta$ is
a surjective map. Also, we have:
$(\rho_r \circ \xi_r)(x) = \rho_r(\xi_r(x)) = \rho_r(x_r) = x_r + I_{p^r, n} =
\Xi_r((x_r + I_{p^r,n})_{r\in{\Bbb N}}) = \Xi_r (\rho_r(x_r)) =
(\Xi_r \circ \varprojlim\rho_r)(x)$. Hence the proposition follows.
\end{proof}

Proposition \ref{prodense} tells us that, although $ {\rm Y}_{\infty ,n}(u)$ may not 
arise as a quotient of $\varprojlim{\Bbb C}{\mathcal F}_{p^r,n}$, yet it does contain a dense 
 quotient. This means that, if we find  a set of topological generators  for 
$\rho ( \varprojlim{\Bbb C}{\mathcal F}_{p^r,n})$
 we will have a set of topological generators  for $ {\rm Y}_{\infty ,n}(u)$. 
In order to do that, we define first certain idempotents $e_{i,j}$ in
$\varprojlim {\Bbb C}{\mathcal F}_{p^r,n}$ that play analogous role to  the
idempontent $e_{p^r,i,j}$. According to (3) in Proposition \ref{proeri} we can define
the following elements:
\begin{equation}\label{eij}
e_{i,j}:= (e_{p,i,j}, e_{p^2,i,j}, \ldots )
\in\varprojlim {\Bbb C}C_r^n \subseteq \varprojlim {\Bbb C}{\mathcal F}_{p^r,n}
\end{equation}
where $i,j \in\{1, \ldots n-1\}$ and $i\not=j$. For  $j=i+1$ we shall denote:
$$
e_i := e_{i,i+1}.
$$
Notice that $e_{i,j} = e_{j,i}$, According to Remark 5, $e_{p^r,i,j}$ is also 
an element in $ {\rm Y}_{p^r, n}(u)$. So (\ref{eij}) defines an element in 
 ${\rm Y}_{\infty, n}(u)$ (with same notation) and we have from the diagram below: 
$$
\begin{CD}
\varprojlim {\Bbb C}C_r^n @>>> \varprojlim {\Bbb C}{\mathcal F}_{p^r, n}  @>\rho >>{\rm
Y}_{\infty  ,n}(u)  \\
 @VVV         @V\xi_r VV     @V \Xi_rVV  \\
{\Bbb C}C_r^n @>>>  {\Bbb C}{\mathcal F}_{p^r, n}  @>\rho_r >>{\rm Y}_{p^r ,n}(u) @>>> 1 \\
\mbox{}
\end{CD}
$$
$$
(\Xi_r \circ \rho ) (e_{i,j}) = (\rho_r\circ \xi_r )(e_{i,j}) = e_{p^r,{i,j}}, 
\quad (\text{for all} \,\, r).
$$

\begin{lem}\label{propeprij}
For any $i, j$ with $i\not= j$, the elements $e_{i,j}\in \varprojlim {\Bbb C}{\mathcal
F}_{p^r,n}$ are idempotents.
\end{lem}
\begin{proof}
The multiplication in $\varprojlim {\Bbb C}{\mathcal F}_{p^r,n}$
is defined componentwise, so the proof follows directly from Lemma \ref{ideeri}.
\end{proof}

\begin{lem}
In  $\varprojlim {\Bbb C}{\mathcal F}_{p^r,n}$, we have:
$$
\sigma_i^2 = 1 + (u-1)e_i (1-\sigma_i) \quad {\rm mod}\, (\varprojlim I_{p^r, n}).
$$
\end{lem}
\begin{proof}
We must prove that $\sigma_i^2 - ( 1 + (u-1)e_i (1-\sigma_i)) \in\varprojlim I_{p^r, n}$. Recall that
$\sigma_i$ is the constant sequence $(\sigma_i, \sigma_i, \ldots )$, hence $\sigma_i^2$ is the constant
sequence $(\sigma_i^2, \sigma_i^2, \ldots )$. Also, the $r$th component of the element 
$ 1 + (u-1)e_i (1-\sigma_i)\in \varprojlim {\Bbb C}{\mathcal F}_{p^r,n}$  is 
$1 + (u-1)e_{p^r,i} (1-\sigma_i)\in{\Bbb C}{\mathcal F}_{p^r,n}$. Therefore, 
the element $\sigma_i^2 - ( 1 + (u-1)e_i (1-\sigma_i))$ is the sequence
$(\epsilon_{1,r}, \epsilon_{2,r},\ldots)$, and $\epsilon_{i,r}\in I_{p^r,n}$. Hence the lemma
follows.
\end{proof}

\begin{prop}
Setting 
$\boldsymbol{\epsilon}_i : =  
{\sigma}_i^2  -1 - (u-1)e_{i} +(u-1)e_{i}{\sigma_i}\in \varprojlim {\Bbb C}{\mathcal
F}_{p^r,n}$,   we have:
$$
\varprojlim I_{p^r, n} =
\langle \boldsymbol{\epsilon}_i \,; \, i=1,\ldots , n-1  \rangle .
$$
\end{prop}
\begin{proof}
Recall  that $\boldsymbol{\epsilon}_i = (\epsilon_{r,i})_{r\in {\Bbb N}}$. Now, for any $i$ and for any 
$x = (x_r)$, $y= (y_r) \in \varprojlim{\Bbb C}{\mathcal F}_{p^r,n}$ we have that
 $x\boldsymbol{\epsilon}_i y = (x_r\epsilon_{r,i}\, y_r)$.  Furthermore  $\phi_s^r(x_r\epsilon_{r,i}y_r)
 = \phi_s^r(x_r)\epsilon_{s,i}\phi_s^r(y_r)\in I_{p^s,n}.$ Thus, $x\boldsymbol{\epsilon}_i y$ belongs to
 $\varprojlim I_{p^r, n}$ for all $i$. Hence, the ideal generated by the $\boldsymbol{\epsilon}_i$'s is
contained in $\varprojlim I_{p^r, n}$. Let now $w = (w_r)_{r\in {\Bbb N}}\in \varprojlim I_{p^r, n}$.
Then  $w_r = \sum_{i} y_{r,i}\epsilon_{r,i}z_{r,i}$, where  $y_{r,i}$,
$z_{r,i}\in {\Bbb C}{\mathcal F}_{p^r,n}.$ Thus, we can write:
$$
w = \sum_i(y_{r,i})_r (\epsilon_{r,i})_r(z_{r,i})_r \in \varprojlim I_{p^r, n}.
$$
As $(y_{r,i})_r$, $(z_{r,i})_r \in \varprojlim{\Bbb C}{\mathcal F}_{p^r,n}$ we obtain $w\in 
\langle \boldsymbol{\epsilon}_i  \, ; \, i=1,\ldots, n-1\rangle$.
\end{proof}

Recall that, according to our inverse system, the element $\sigma_i\in B_n$
corresponds to the constant sequence $(g_i, g_i, \ldots )$ in ${\rm Y}_{\infty ,n}(u)$.  We denote  this
 sequence by $g_i$. Similarly, the braid $\sigma_i^{-1}\in B_n$ corresponds to the
 constant sequence $(g_i^{-1}, g_i^{-1}, \ldots )$ in ${\rm Y}_{\infty ,n}(u)$ and 
it shall be denoted by $g_i^{-1}$. Thus, in $\rho ( \varprojlim{\Bbb C}{\mathcal F}_{p^r,n})
\subseteq {\rm Y}_{\infty,n}(u)$ the following quadratic relations holds:
$$
g_i^2 = 1 + (u-1)e_i
(1- g_i)\quad (i=1, \ldots ,n-1 ).
$$
We define now  ${\bold t}_i:= \rho({\bold t}_i)$  and $ e_i:= \rho(e_i)$. Then, from  Theorem \ref{yrnr}
and Proposition \ref{prodense}, we deduce  the following theorem.

\begin{thm}\label{toppre}
 $\{1, {\bold t}_1, g_1,\ldots , g_{n-1}\}$ is a set of topological generators of ${\rm
Y}_{\infty,n}(u)$. Moreover, these elements satisfy the following relations:
 \begin{enumerate}
\item Braid relations among the $g_i$'s
\item ${\bold t}_1 g_i =g_i {\bold t}_1$, for $i\geq 2$
\item ${\bold t}_1 g_1{\bold t}_1 g_1^{-1} =
  g_1{\bold t}_1 g_1^{-1}{\bold t}_1$
\item
$g_ig_{i-1}\cdots
 g_1{\bold t}_1g_1^{-1}\cdots g_{i-1}^{-1})
 g_i^{-1}
= g_i^{-1}(g_{i-1}\cdots 
g_1{\bold t}_1g_1^{-1}\cdots
g_{i-1}^{-1})g_i$
\item $g_i^2 = 1 + (u-1)e_i (1-g_i)$,
$(i=1,\ldots , n-1).$
\end{enumerate}
\end{thm}
Moreover, as in  Proposition \ref{peri},  we  can  prove analogous commutation relations for $e_i$.
More precisely:

\begin{prop}\label{pei}
In ${\rm Y}_{\infty,n}(u)$ we have :
\begin{enumerate}
\item $g_i^{\pm 1} e_j=
e_jg_i^{\pm 1}$, for $j\not=i-1, i+1$
\item $g_i^{\pm 1} e_j = e_{ij} g_i^{\pm 1}$, 
for $\vert i-j\vert =1$
\item $e_j g_i^{\pm 1} = g_i^{\pm 1} e_{ij} $, 
for $\vert i-j\vert =1$.
\end{enumerate}
\end{prop}
\begin{proof}
The proofs  follow directly from Lemma \ref{propeprij} and Proposition
\ref{peri}.
\end{proof}

\begin{rem}\rm 
 It is worth  observing  that  ${\rm Y}_{\infty, n}(u)$ can be 
 regarded as a topological deformation of a quotient of the group algebra 
 ${\Bbb C}{\mathcal F}_n$, see Theorem \ref{toppre}. 
Roughly, the algebra ${\rm Y}_{\infty, n}(u)$ can 
 be described in terms of topological generators, in the sense of Definition \ref{generators}, 
and the same  relations as the algebra  ${\rm Y}_{d , n}(u)$ but where the relations
$h_i^{d}= 1$ do not hold. Consequently, ${\rm Y}_{\infty, n}(u)$
 has a set a toplogical generators, which look like  the
canonical generators  of the framed braid group ${\mathcal F}_n$ (recall Proposition
\ref{prefn1}), but with the addition of the quadratic relation. 
\end{rem}

\section{Topological Markov traces}

In \cite{ju} the first author constructed linear Markov traces on the Y--H algebras. The aim of this
section is to extend these traces to a Markov trace on the algebra ${\rm Y}_{\infty ,n}(u).$

\subsection{}

The natural inclusions of the framed braid groups: ${\mathcal F}_{i} \subset
{\mathcal F}_{i+1}$ induce natural inclusions of the modular framed
braid groups: ${\mathcal F}_{k,i} \subset {\mathcal F}_{k,i+1},$ which in turn
induce the tower of algebras:
$$
{\rm Y}_{k,0} (u) := {\Bbb C} \, \subset {\rm Y}_{k,1} (u) \subset {\rm Y}_{k,2} (u) \subset \ldots
$$
We then have the following.

\begin{thm}[cf. Theorem 12 in \cite{ju}]\label{tyk}
For indeterminates $z$, $x_1$, $\ldots, x_{k-1}$ there exists a unique linear Markov trace
$$
 {\rm tr}: \ {\rm Y}_{k, n+1}(u)  \longrightarrow   {\Bbb C}[z, x_1, \ldots, x_{k-1}]
$$
 defined inductively by the following rules:
$$
\begin{array}{rcll}
{\rm tr}(ab) & = & {\rm tr}(ba)  \qquad &  \\
{\rm tr}(1) & = & 1 & \\
{\rm tr}(ag_n b) & = & z\, {\rm tr}(ab) \qquad &  \\
{\rm tr}(ah_{n+1}^m b) & = & x_m {\rm tr}(ab)\qquad  & (  m = 1, \ldots  , k-1),
\end{array}
$$
where $a,b \in {\rm Y}_{r,n}(u)$.
\end{thm}

\subsection{\it The $p$-adic Markov trace}

Let ${\rm X}_r= \{z, x_1, x_2, \ldots ,x_{p^r -1}\}$ be a set of indeterminates. For $r>s$ we 
define the `connecting' ring epimorphism: 
$$
\delta_s^r :{\Bbb C}[{\rm X}_r] \longrightarrow
{\Bbb C}[{\rm X}_s]
$$ 
as the composition
$$
{\Bbb C}[{\rm X}_r] \longrightarrow \frac{{\Bbb C}[{\rm X}_r]}{{\mathcal I}^r_s} \simeq {\Bbb
C}[{\rm X}_s],
$$
where ${\mathcal I}_s^r$ is the ideal of ${\Bbb C}[{\rm X}_r]$ generated by the polynomials:
$$
\begin{array}{clll}
x_i -x_j & \text{for} & i \equiv j & \, {\rm mod} (p^s) \\
x_i - 1  & \text{for} & i \equiv 0 & \, {\rm mod} (p^s),
\end{array}
$$
for any $i,j \in \{1, \ldots , p^r -1\}$.  More precisely, the map $\delta_s^r$ acts on the indeterminates as
follows:
$$
\delta_s^r(x_i) =
\left\{
\begin{array}{ll}
x_i, & \text{for}\, i= 1, \ldots , p^s-1\\
1, & \text{for}\, i\equiv0 \, ({\rm mod}\, p^s)\\
x_j, & \text {for}\,i\equiv j \, ({\rm mod}\, p^s).
\end{array}\right.
$$

It is then  a  routine to  prove the following lemma.
\begin{lem}
The family $\left( {\Bbb C}[X_r], \delta_s^r \,;\, r\in {\Bbb N}\right)$ is an
inverse system of polynomial rings.
\end{lem}

Now let  $\tau_r={\rm tr}$,  the trace on the Y--H algebra  ${\rm
Y}_{p^r,n}(u),$  as defined in Theorem~\ref{tyk} above, taking
 values in ${\Bbb C}[{\rm X}_r]$. Then,  $\tau_r(g_n) = z$ and 
 $\tau_r(h_n^m)= x_m$,  where $m = 1, \ldots , p^r-1$.
 We then  obtain the following commutative diagram:
$$
\begin{CD}
 {\rm Y}_{p^r,n}(u) @>\pi_s^r>> {\rm Y}_{p^s,n} (u)\\
   @V \tau_r VV         @V\tau_sVV    \\
{\Bbb C}[{\rm X}_r] @> \delta_s^r >> {\Bbb C}[{\rm X}_s]
\end{CD}
$$
It follows that the family $E=\{\tau_r \}_{r\in {\Bbb N}}$ of homomorphisms of inverse systems, 
$$
E: \left({\rm Y}_{p^r,n}(u), \pi_s^r\right) \longrightarrow 
\left({\Bbb C}[{\rm X}_r], \delta_s^r\right)
$$
yielding a unique ring homomorphism  $\varprojlim \tau_r$. Thus, we have the following theorem.

\begin{thm}\label{ptrace}
There  exists a unique $p$-adic linear Markov trace defined  as
 $$
\tau := \varprojlim \tau_r :{\rm Y}_{\infty ,n+1}(u)
\longrightarrow
\varprojlim {\Bbb C}[{\rm X}_r]$$
Furthermore 
$$
\begin{array}{rcl}
\tau(ab)  & = & \tau(ba)     \\
\tau (1) & = & 1  \\
\tau(ag_nb) & = & (z)_r \tau(ab)    \\
\tau (a{\bf t}_{n+1}^mb) & = & (x_m)_r \tau(ab)  
\end{array}
$$
for any $a,b \in {\rm Y}_{\infty ,n}(u)$ and $m\in {\Bbb Z}$.
\end{thm}
\begin{proof}
By  construction, the trace $\tau$ satisfies all the properties in the statement.
\end{proof}

\subsection{}

In this subsection we shall give some computations of the trace $\tau$. 
For  a $p$--adic integer $ \underleftarrow{a} = (a_1, a_2, \ldots) \not = 0$, we shall denote
$$
x_{\underleftarrow{a}}= (x_{a_1}, x_{a_2}, \ldots ) \in  \varprojlim {\Bbb C}[{\rm X}_r].
$$
We call  $x_{\underleftarrow{a}}$ a $p$--{\it adic indeterminate}. In the case 
$\underleftarrow{a} = (k, k, \ldots)\in {\Bbb Z}\subset {\Bbb Z}_p$, we denote 
 $x_{\underleftarrow{a}}$ by $x_{k}$ and we say that $x_k$ is a {\it constant indeterminate}. Finally, we 
set  $z:=(z,z, \ldots )$ and we also make the convention $x_0:=1$.

\smallbreak
\noindent $\bullet$ 
For ${\bf t}\in {\rm Y}_{\infty ,1}(u)$ we have:
$$
\tau ({\bf t}) = \varprojlim \tau_r ({\bf t})= (\tau_1(t_1), \tau_2(t_2), \ldots )=  (x_1,
x_1,\ldots ) = x_1.
$$
\noindent $\bullet$ 
We shall now compute  the trace of an element $w={\bf t}^{\underleftarrow{a}}\in {\Bbb Z}_p
$. We can write  $w$ as $w = (t_1^{a_1}, t_2^{a_2}, \ldots)$. Then:
$$
\tau ({\bf t}^{\underleftarrow{a}}) = (\tau_1(t_1^{a_1}), \tau_2(t_2^{a_2}), \ldots )
= (x_{a_1}, x_{a_2}, \ldots ) = x_{\underleftarrow{a}}.
$$
Further, since $ {\bf t}^{\underleftarrow{a}} =  \lim_{k} {\bf t}^{a_k}$, we have the following
aproximation:
$$
\tau({\bf t}^{\underleftarrow{a}}) = \lim_{k} \tau ({\bf t}^{a_k})  = \lim_{k} x_{a_k}.
$$
\noindent $\bullet$  For  $m_i\in {\Bbb Z}$  we have:
$$
\tau ({\bf t}_1^{m_1}\cdots {\bf t}_n^{m_n} )  = \tau ({\bf t}_1^{m_1})\cdots
\tau({\bf t}_n^{m_n} ) = x_{m_1}\cdots x_{m_n}. 
$$
\noindent$\bullet$ 
In general, let $w\stackrel{!}{=}  {\bf t}_1^{\underleftarrow{a_1}} {\bf t}_2^{\underleftarrow{a_2}}\cdots 
 {\bf t}_n^{\underleftarrow{a_n}}$, where $\underleftarrow{a_i} = (a_{ri})_r$. Then we have:
$$
\tau(w) = \tau ({\bf t}_1^{\underleftarrow{a_1}}    {\bf t}_2^{\underleftarrow{a_2}}\cdots 
 {\bf t}_n^{\underleftarrow{a_n}} )= 
\tau ({\bf t}_1^{\underleftarrow{a_1}})\cdots \tau({\bf t}_n^{\underleftarrow{a_n}}) 
= x_{\underleftarrow{a_1}}\cdots
x_{\underleftarrow{a_n}}.
$$
Indeed, $w =(t_{r,1}^{a_{r1}} \, t_{r,2}^{a_{r2}} \, \cdots \, t_{r,n}^{a_{rn}})_r \in {\rm Y}_{\infty ,n}(u)
$. Hence:  
\begin{eqnarray*}
\tau (w) & = & ( \tau_1(t_{1,1}^{a_{11}}\cdots t_{1,n}^{a_{1n}}) ,
\tau_2( t_{2,1}^{a_{21}}\cdots t_{2,n}^{a_{2n}}),\ldots) \\
 &  = & 
(x_{a_{11}}\cdots x_{a_{1n}},x_{a_{21}}\cdots x_{a_{2n}}, \ldots )\\
 & = & 
(x_{a_{11}},x_{a_{21}},\ldots) \cdots (x_{a_{1n}},x_{a_{2n}},\ldots)  \\
 & = & 
 x_{\underleftarrow{a_1}}\cdots x_{\underleftarrow{a_n} }  \\
 & = & 
 \tau ({\bf t}_1^{\underleftarrow{a_1}})\cdots \tau({\bf
t}_n^{\underleftarrow{a_n}}).
\end{eqnarray*}
Further, for $\tau (w)$  we also have the approximation:
$$
\tau({\bf t}_1^{\underleftarrow{a_1}}    {\bf t}_2^{\underleftarrow{a_2}}\cdots 
 {\bf t}_n^{\underleftarrow{a_n}}) = 
\lim_{k} \tau ({\bf t}_1^{a_{k1}} \cdots  {\bf t}_n^{a_{kn}} ) = 
\lim_{k} \left(x_{a_{k1}} \cdots  x_{a_{kn}} \right). 
$$

In order to compute the trace of $g_i^2$, we first compute the trace 
 at the level of the classical Y--K algebras. Indeed, recall that 
$g_i^2 = 1 + (u-1)e_{p^r, i} - (u-1)e_{p^r, i}g_i$. 

\smallbreak
\noindent $\bullet$ We shall first compute  $\tau_r (e_{p^r, i}g_i)$ and $\tau_r (e_{p^r, i})$. Indeed, we
have:
\begin{eqnarray*}
\tau_r (e_{p^r, i}g_i) & = & \frac{1}{p^r}\sum_{m=0}^{p^r-1}\tau_r (t_{r,i}^mt_{r, i+1}^{-m}g_i) \\
& = & 
\frac{1}{p^r}\sum_{m=0}^{p^r-1}\tau_r (t_{r,i}^mg_it_{r, i}^{-m})\quad (\text{from} \, (\ref{tiis}))\\
& = & 
\frac{1}{p^r}\sum_{m=0}^{p^r-1}z \, \tau_r (t_{r,i}^m t_{r, i}^{-m}) = z.
\end{eqnarray*}
Thus 
$$
\tau_r (e_{p^r, i}g_i) = z. 
$$
Moreover,
$$
\tau_r(e_{p^r, i}) 
= 
\tau_r\left(\frac{1}{p^r}\sum_{m=0}^{p^r-1}t_{r,i}^mt_{r, i+1}^{-m} \right)
=
\frac{1}{p^r}\sum_{m=0}^{p^r - 1}x_{m}x_{-m}.
$$

\noindent $\bullet$ Then, for $\tau_r (g_i^2)$ we have: 

\begin{eqnarray*}
\tau_r (g_i^2) & = &  1 + (u-1)\tau_r (e_{p^r, i})- (u-1)\tau_r (e_{p^r, i}g_i)  \\
 & \Updownarrow  & \\
\tau_r(g_i^2) & = & 1 - (u-1)z + (u-1)\frac{1}{p^r}\sum_{m=0}^{p^r - 1}x_{m}x_{-m}.
\end{eqnarray*}

In order to compute now the trace of $g_i^2= (g_i^2, g_i^2,\ldots )\in  
{\rm Y}_{\infty ,n}$, we need to compute first $\tau(e_i)$ and $\tau(e_ig_i)$. 

\smallbreak
\noindent $\bullet$ We have $e_i = (e_{p^r,i})_r$, so: 
$$
\tau(e_i) = 
\left(
\frac{1}{p}\sum_{m=0}^{p - 1} x_{m}x_{-m} \, , 
\frac{1}{p^2}\sum_{m=0}^{p^2 - 1} x_{m}x_{-m} \, , \ldots
\right).
$$
Also, a  direct computation gives: 
$$
\tau (e_ig_i) = z.
$$
\noindent $\bullet$ Recall that 
$$
g_i^2 = 1 + (u-1)e_i - (u-1)e_ig_i. 
$$
From the above, we obtain for $\tau (g_i^2 )$: 
$$
\tau (g_i^2 )= 1- (u-1)z + (u-1)\left(\frac{1}{p^r}\sum_{m=0}^{p^r -1} x_{m}x_{-m} \right)_r,
$$
where $1:= (1, 1, \ldots )$.

\smallbreak
\noindent $\bullet$ Finally, we shall give approximations for $\tau(e_i)$ and $\tau (g_i^2 )$. 
Let us define
$$
z_{r,k,i}
:=\frac{1}{p^k}\sum_{m=0}^{p^k}t_{r,i}^mt_{r, i+1}^{-m} 
\stackrel{!}{=} 
\frac{1}{p^k}\sum_{m=0}^{p^k}{\bf t}_i^m{\bf t}_{ i+1}^{-m}
$$
Notice that $z_{r,r,i} = e_{{p^r},i}$. Considering now the sequence $(z_{r,k,i })_k$ of constant sequences, we
have:
$$
e_i = \lim_{k} z_{r, k,i}. 
$$
Further, 
$$
\tau (z_{r, k,i}) = \frac{1}{p^k}\sum_{m=0}^{p^k-1}\tau({\bf t}_i^m {\bf t}_{i+1}^{-m}) =
\frac{1}{p^k}\sum_{m=0}^{p^k-1}x_m x_{-m} \in \varprojlim {\Bbb C}[{\rm X}_r].
$$
Therefore,
$$
\tau(e_i) = \lim_{k}\left(\frac{1}{p^k}\sum_{m=0}^{p^k-1}x_m x_{-m}\right).
$$
A direct computation also gives:
$$
\tau(e_ig_i) 
= \tau \left(\lim_{k} \left(z_{r, k,i}g_i\right)_k\right) 
=
\lim_{k} (z)_k = z.
$$
So, we obtain
$$
\tau (g_i^2) = 1 - (u-1)z + (u-1) \lim_{k}\left(\tau(z_{r,k,i})\right).
$$

\begin{rem}\rm 
As already noted in the Introduction, in a sequel paper we combine the $p$-adic trace $\tau$ of
Theorem~\ref{ptrace}, as well as the Markov traces of
Theorem~\ref{tyk},  with  the Markov  equivalence for $p$-adic framed braids to construct invariants of
oriented $p$-adic framed links. We hope that this new concept of  $p$-adic framed braids and  $p$-adic framed
links that we propose, as well as our $p$-adic framing invariant, will be useful for constructing
$3$--manifold invariants using the theory of braids and the Markov-type equivalence given in \cite{ks}. 
\end{rem}




\begin{thebibliography}{50}


\bibitem{bou} N. Bourbaki, {\em Elements of mathematics. Algebra}. Paris: Herman.
Chapitre 2 (1962).


\bibitem{fr} R. Fenn, C.P. Rourke, {\em On Kirby's calculus of links}, 
Topology {\bf 18}, 1--15 (1979).




\bibitem{jo} V.F.R. Jones, {\em Hecke algebra representations of braid groups and link
polynomials}, Ann. Math. {\bf 126}, 335--388 (1987).

\bibitem{ju} J. Juyumaya, {\em Markov trace on the Yokonuma-Hecke algebra},
J. Knot Theory and its Ramifications {\bf 13}, 25--39 (2004).

\bibitem{ki} R. Kirby, {\em A Calculus for Framed Links in $S^3$}, 
Inventiones Math. {\bf 45}, 35--56 (1978).

\bibitem{ks} K.H. Ko, L. Smolinsky, {\em The framed braid group and $3$--manifolds},
Proceedings of the AMS, {\bf 115}, No. 2, 541--551 (1992).

\bibitem{la} S. Lambropoulou, {\em Knot theory related to generalized and cyclotomic
Hecke algebras of type {\rm B}}, J. Knot Theory and its Ramifications {\bf 8},
No. 5, 621--658 (1999).

\bibitem{li} W.B.R. Lickorish, {\em A representation of orientable combinatorial 3-manifolds}, 
Annals of Mathematics {\bf 76}, No. 3, 531--540 (1962).

\bibitem{riza} L. Ribes and P. Zalesskii, {\em Profinite Groups}, A Ser. Mod. Sur. Math. 40,
Springer (2000).

\bibitem{ro} A.M. Roberts, {\em A Course in $p$-adic Analysis}, Grad. Texts in Math. 198,
Springer (2000).

\bibitem{th} N. Thiem {\em Unipotent Hecke algebras}, 
Journal of Algebra, {\bf 284}, 559--577 (2005).

\bibitem{wa} A.D. Wallace, {\em Modifications and cobounding manifolds}, 
Can. J. Math. {\bf 12} , pp. 503--528 (1960).

\bibitem{wi} Wilson, {\em Profinite Groups}, London Math. Soc. Mono., New Series 19,
Oxford Sc. Publ. (1998).

\bibitem{yo} T. Yokonuma, {\em Sur la structure des anneaux de Hecke d'un
groupe de Chevallley fini,} C.R. Acad. Sc. Paris, {\bf 264},  344--347 (1967).

\end{thebibliography}
\end{document}